\documentclass[journal]{IEEEtran}

\ifCLASSINFOpdf
\else
   \usepackage[dvips]{graphicx}
\fi
\usepackage{url}

\usepackage{amsmath, amssymb, amsthm, float, graphicx, subfig, array, mathrsfs, xr}

\usepackage[ruled,vlined]{algorithm2e}
\SetKwInput{KwInput}{input}
\SetAlFnt{\small}

\newtheorem{theorem}{Proposition}
\newtheorem{lemma}{Lemma}
\newtheorem{remark}{Remark}

\begin{document}

\title{Per-antenna power constraints: constructing Pareto-optimal precoders with cubic complexity under non-negligible noise conditions}

\author{Sergey Petrov, Samson Lasaulce, Merouane Debbah
\thanks{This work is sponsored by the Khalifa University-TII Educational Agreement '6G Chair on Native AI'.}
\thanks{Sergey Petrov is with Khalifa University (e-mail: spetrov.msk@gmail.com, sergey.petrov@ku.ac.ae); Samson Lasaulce and Merouane Debbah are with Khalifa University.}}

\markboth{Journal of \LaTeX\ Class Files, Vol. 14, No. 8, August 2015}
{Shell \MakeLowercase{\textit{et al.}}: Bare Demo of IEEEtran.cls for IEEE Journals}
\maketitle

\begin{abstract}
Precoding matrix construction is a key element of the wireless signal processing using the multiple-input and multiple-output model. It is established that the problem of global throughput optimization under per-antenna power constraints belongs, in general, to the class of monotonic optimization problems, and is unsolvable in real-time. The most widely used real-time baseline is the suboptimal solution of Zero-Forcing, which achieves a cubic complexity by discarding the background noise coefficients. This baseline, however, is not readily adapted to per-antenna power constraints, and performs poorly if background noise coefficients are not negligible. In this paper, we are going to present a computational algorithm which constructs a precoder that is SINR multiobjective {\it Pareto-optimal} under per-antenna power constraints - with a complexity that differs from that of Zero-Forcing only by a constant factor. The algorithm has a set of input parameters, changing which skews the importance of particular user throughputs: these parameters make up an efficient parameterization of the entire Pareto boundary.
\end{abstract}


\IEEEpeerreviewmaketitle

\section{Introduction}
\subsection{General MIMO precoding optimization problem}

\IEEEPARstart{T}{he} transmission of a signal between a transmitter with $m_{tx}$ antennas and a $m_{ue}$ receiving antennas (or spatial streams), according to the wireless multiple-input multiple-output (MIMO)\footnote{The developments of the paper can be applied to $m_{ue}$ receivers with one antenna or one receiver with $m_{ue}$ antennas with no mathematical difference.} model, can be described with an equation
\begin{align*}
y & = H^* x + n, H \in \mathbb{C}^{m_{tx} \times m_{ue}} \\
x & = P s, P \in \mathbb{C}^{m_{tx} \times m_{ue}},
\end{align*}
where $s$ is the complex-valued signal intended to be sent, $x$ is a preprocessed signal, $n$ is the noise, and $y$ is the received signal. The so-called 'channel' matrix $H$ is determined by physical properties of the transmission and is predefined\footnote{$H^*$ denotes a conjugate transpose of the matrix $H$.}, while the 'precoding' matrix $P$ is subject to be selected. Through this paper, unless directly specified otherwise, for any capital latin letter matrix notation $A$ the notations $a_j$ and $a_{ij}$ will be used for column vectors and elements of $A$, respectively. The quality of a precoding matrix $P$ is then commonly measured through 'Signal-to-Interference-and-Noise Ratio, SINR' objective values as:
\begin{align}
	G & := H^* P, G \in \mathbb{C}^{m_{ue} \times m_{ue}}; \nonumber \\
	{\cal S}_k & := \lvert g_{kk} \rvert^2, {\cal I}_k := \sum \limits_{j \neq k}^{m_{ue}} |g_{kj}|^2;  \nonumber \\
	{\cal N}_k & := \omega_k^2, {\cal W}_k := {\cal I}_k + {\cal N}_k; \nonumber \\
	\mathrm{SINR}_k(P) & := \frac{{\cal S}_k}{{\cal W}_k} = \frac{|h_k^* p_k|^2}{\omega_k^2 + \sum \limits_{j = 1, j \neq k}^{m_{ue}} |h_k^* p_j|^2}; \label{eq:sinr}
\end{align}
The precoder $P$ must satisfy a set of practical 'power' constraints, which come in the following 'harder' and 'simpler' variants:
\begin{align}
	\| e_i^* P \|_2^2 & \leq \beta_i;
	\mbox{ (Per-antenna, practical)} \label{eq:powerrow} \\
	\| P \|_F^2 & \leq \sum \limits_{i = 1}^{m_{tx}} \beta_i,
	\mbox{ (Frobenius, simplified)} \label{eq:powerfro}
\end{align}
where $e_i$ is a canonical unit vector with only one nonzero element equal to $1$ at position $i$ and $\beta > 0$ are power limit values. We will call a collection of objective values (\ref{eq:sinr}) {\it Pareto-optimal} if one of those cannot be increased without reducing the values of others or violating the selected power constraint. We will call a precoding matrix $P$ Pareto-optimal if the corresponding combination of $\mathrm{SINR}_k(P)$ is Pareto-optimal. In this paper, we will focus on the more complicated per-row power constraint type (\ref{eq:powerrow}), which is actual when each antenna is fed by its own dedicated power amplifier.

Whenever we shall need to select {\it one} particular Pareto-optimal point, we will attempt to additionally optimize to the total/mean throughput 
\begin{equation}
\mathscr{T}(P) := \sum \limits_{k = 1}^{m_{ue}} \log(1 + \mathrm{SINR}_k(P)), \label{eq:func}
\end{equation}
which we will refer to as 'global optimization'. We will subject ourselves to the following limitations:
\begin{itemize}
	\item Perfect knowledge of the channel matrix $H$ is assumed. 
	The 'non-negligible noise' of the title refers to the coefficients $\omega_k$, but not the channel estimation quality.
	\item The number of transmitter antennas dominates the number of receiver antennas: $m_{tx} > m_{ue}$. It is important from the complexity point of view: generally speaking, a complexity of $O(m_{tx}^2 m_{ue})$ will be seen as strictly worse than $O(m_{tx} m_{ue}^2)$.
\end{itemize}

If the transmission noise variance $\omega^2$ could be ignored, the optimal way to select $P$ would be such that $H^* P$ is a diagonal matrix. Such a selection
\begin{equation*}
\label{eq:zf}
P_{\mathrm{ZF}} := H (H^* H)^{-1}
\end{equation*} 
is known as a 'Zero-Forcing' (ZF) \cite{zf} solution and often serves as a baseline.  Since wireless applications commonly have real-time requirements, the low finite complexity $O(m_{tx} m_{ue}^2)$ makes zero-forcing an appealing strategy. A lot of research has been devoted to reducing its complexity even further: works \cite{zhang1, zhang2} utilize Neumann series, fast Fourier transform and block-Toeplitz algorithms to construct fast approximations for $H(H^* H)^{-1}$, while papers \cite{mitya, kavi} introduce blocks and work on pushing a low-complexity block diagonal zero-forcing towards the performance of full dense $P_{\mathrm{ZF}}$.

For the zero-forcing baseline, the power constraints are commonly fulfilled by post-weighting the columns of $P_{\mathrm{ZF}}$, so as not to disrupt the property of $H^* P_{\mathrm{ZF}}$ being diagonal. This column weighting coefficient optimization for the case (\ref{eq:powerfro}) is known to be solved via an efficient water-filling algorithm \cite{water}. Weighting gets more complicated with per-antenna case, but is also well-studied: work \cite{bocca} utilizes interior-point methods to optimize the scaling coefficients for (\ref{eq:func}), gives analytical expressions for $m_{ue} = 2$ and shows than, in most cases, water-filling can still be viable when applied to one of the per-antenna constraints while discarding the others. Work \cite{vu} considers the power limits $\beta$ as additional variables, and aims at zero-forcing weighting for the optimal rate between ($\ref{eq:func}$) and power used instead.

Optimizing column weights for $P_{\mathrm{ZF}}$ can have undesirable additional computational complexity. In practical applications with per-antenna constraints (\ref{eq:powerrow}), $P_{\mathrm{ZF}}$ is commonly simply scaled by a matrix-wide constant 
\begin{equation*}
P_{\mathrm{ZF}} \xrightarrow{\text{scaling}} \min \limits_i \frac{\sqrt{\beta_i}}{\|e^*_i P_{\mathrm{ZF}} \|_2} P_{\mathrm{ZF}}.
\end{equation*}
The work \cite{lee} gives some bounds on SINR distributions and an asymptotic ($m_{tx}, m_{ue} \rightarrow \infty$) performance study of such a scaling approach for random Gaussian channels. The approach, however, fulfills only one per-antenna constraint with equality, hence 'full power' is not used. That raises one of the questions to be studied in this paper: we claim that full power usage is a {\it necessity} for all Pareto-optimal precoders under mild assumptions on the channel. 

Consequently, even if column weights for $P_{\mathrm{ZF}}$ are selected {\it optimally} with respect to (\ref{eq:func}), the resulting precoding matrix is still {\it suboptimal} for (\ref{eq:func}): most importantly, because it is restricted to the column subspace of $H$. Informally, fulfilling {\it all} the per-antenna power constraints with equality with a precoder restricted to the channel column subspace is impossible in general: it corresponds to $m_{tx}$ restrictions with $O(m_{ue}^2)$ precoding matrix subspace coordinate variables. Formally, for Frobenius \cite{wiesel}\cite{bjornsson} or per-user \cite{jors} power constraints, it is shown that using the channel column space for precoding is {\it sufficient}. It is {\it not} the case for per-antenna power constraints (\ref{eq:powerrow}): the work \cite{wiesel} presents some generalization of zero-forcing that escapes the column space of $H$, but keeps $H^*P$ diagonal, and brings a positive throughput gain under constraints (\ref{eq:powerrow}).

The global optimizer of (\ref{eq:func}) under per-antenna constraints (\ref{eq:powerrow}), in general, fulfills neither the $H^* P$-diagonality, nor the column subspace $span(H) = span(P)$ assumptions. Without those, the problem loses convexity \cite{nonconv1, nonconv2}.
It can still be solved using a combination of a monotonic optimization method, such as 'Branch-Reduce-Bound' \cite{brb}, and a convex solver for the so-called 'Quality-of-Service' (QoS) problem \cite{bjornsson}, but such an optimization is on a different level of complexity: not only it is inapplicable in real time, but it can hardly converge in hour/day time scale for $m_{ue} > 10$. 

In this work, we will consider the less studied scenario when $H^* P$ cannot be limited to a diagonal matrix due to higher noise coefficients $\omega$. We will not be aiming at optimizing (\ref{eq:func}) globally; instead, we will focus on efficiently computing a family of Pareto-optimal precoders under per-antenna power constraints (\ref{eq:powerrow}) in real-time cubic complexity. We will argue that this computation can be used both as a self-sufficient real-time solution and as an auxiliary tool for global optimization.

\subsection{Fixed Quality-of-Service convex problem}
\label{subsec:qos}
The fixed-QoS problem corresponds to establishing whether the solution set of an inequality system for a given $\gamma_1 \dots \gamma_{m_{ue}}$
\begin{equation}
	\label{eq:qos}
	\begin{cases}
		\mathrm{SINR}_1(P) \geq \gamma_1 \\
		\mathrm{SINR}_2(P) \geq \gamma_2 \\
		\vdots \\
		\mathrm{SINR}_{m_{ue}} \geq \gamma_{m_{ue}}
	\end{cases}
\end{equation}
is nonempty under the power constraints (\ref{eq:powerrow}). This problem was studied in depth in \cite{bjornsson}: utilizing $\lvert \alpha \rvert^2 = \alpha \alpha^*, \alpha \in \mathbb{C}$, and recalling (\ref{eq:sinr}), the set of inequalities (\ref{eq:qos}), (\ref{eq:powerrow}) can be rewritten as a convex optimization problem
%\footnote{Through the paper, indeces $j$ and $k$ are swapped sometimes, but there is no way to completely avoid that in some derivations with sum reordering.}
\begin{align}
	\mathrm{minimize} & \mbox{ } 0;  \label{eq:convexbasic} \\
	\mathrm{subject} \mbox{ } \mathrm{to} \mbox{ } \sum \limits_{k = 1}^{m_{rx}} p_k^* e_i e_i^* p_k  - \beta & \leq 0, i = 1 \dots m_{tx}; \nonumber \\
	\gamma_j (\sum \limits_{k \neq j} p_k^* h_j h_j^* p_k + \omega_j^2) - p_j^* h_j h_j^* p_j & \leq 0, j = 1 \dots m_{ue}. \nonumber
\end{align}
By analyzing the corresponding Lagrangian and its derivatives, \cite{bjornsson} (Corollary 3.6) formally proves that any obtainable combination $\gamma_1, \gamma_2, \dots,  \gamma_{m_{ue}}$ of SINR values is necessarily obtainable by a precoder $P$ built with formulas
\begin{align}
	\Psi & := \begin{bmatrix}
		\frac{\mu_1}{\beta_1} & & \\
		& \ddots & \\
		& & \frac{\mu_{m_{tx}}}{\beta_{m_{tx}}} \\
	\end{bmatrix}
	+ \sum \limits_{k = 1}^{m_{ue}} \frac{\lambda_k}{\omega_k^2} h_k h_k^*; \hat p_j := \frac{\Psi^{-1} h_j}{\| \Psi^{-1} h_j \|_2}. \nonumber \\
	\gamma_j & := \frac{\lambda_j}{\omega_j^2} h_j^* (\Psi - \frac{\lambda_j}{\omega_j^2} h_j h_j^*)^{-1} h_j; (S)_{ij} :=
	\begin{cases}
		\lvert h_i^* \hat p_i \rvert^2, i = j \\
		- \gamma_i \lvert h_i^* \hat p_j \rvert^2, i \neq j.
	\end{cases} \nonumber \\
	\kappa & := S^{-1}
	\begin{bmatrix}\gamma_1 \omega_1^2 & \dots & \gamma_{m_{ue}} \omega_{m_{ue}}^2
	\end{bmatrix}^\mathrm{T}; p_j = \sqrt{\kappa_j} \hat p_j. \label{eq:noncomp}
\end{align}
with {\it some} combination of Lagrange multipliers $\mu, \lambda$ with no additional scaling. The formulas (\ref{eq:noncomp}) are, however, numerically complex: $O(m_{tx}^3 m_{ue})$ operations are required if they are applied directly. Also, not every combination $\mu, \lambda$ yields a valuable $P$: selecting the multiplier values at random most commonly results in a violation of power constraints. The Pareto-boundary, assuming smoothness, should be a surface that can be described using $m_{ue} - 1$ parameters, but the Lagrange multipliers $\lambda$ and $\mu$ constitute a set of $m_{tx} + m_{ue}$ parameters. One of those can be discarded using
\begin{equation*}
\sum \limits_{i = 1}^{m_{tx}} \mu_i = \sum \limits_{j = 1}^{m_{ue}} \lambda_j,
\end{equation*}  
which is proved in \cite{bjornsson} by recalling that the optimization functional in (\ref{eq:convexbasic}) is trivial, but $m_{tx}$ 'extra' degrees of freedom remain: $m_{tx}$ additional equations should be imposed if we want (\ref{eq:noncomp})  to become a true Pareto boundary parameterization. 

Our contribution will be organized as follows. First, we prove that {\it full power usage} is mandatory for all Pareto-optimal precoders under some mild realistic assumptions on the channel and an assumption that noise $\omega$ is rather high: each inequality (\ref{eq:powerrow}) must be fulfilled with equality.

Secondly, we will propose a simple numerical scheme to {\it solve the full power constraints} as equalities on $\lambda, \mu$, effectively removing $m_{tx}$ 'extra' degrees of freedom and arriving at an efficient parameterization of a Pareto boundary for the case of per-antenna power constraints.

Finally, we will show that the formulas (\ref{eq:noncomp}) can be equivalently rewritten so that the computational complexity is much more tolerable: we will show that one Pareto-optimal precoding matrix can be computed in constant-times Zero-Forcing complexity of $O(m_{tx} m_{ue}^2)$, making it applicable for real-time applications.
We will then prove the effectiveness of the proposed methods in numerical experiments. 

\section{Full power usage necessity}
\label{sec:power}
In this section, we provide a theorem that guarantees the necessity of the full power usage under certain mild assumptions on the channel $H$ when noise coefficients $\omega$ are not negligible. Let us first annotate the result.
\begin{theorem}{(\it Full power usage necessity.)}
\label{th:power}
	Assume $H \in \mathbb{C}^{m_{tx} \times m_{ue}}$ is a channel matrix such that $m_{tx} > m_{ue} \geq 2$. Assume that a block channel-identity matrix
	\begin{equation*}
		\begin{bmatrix}
			H & \lvert & I_{m_{tx}}
		\end{bmatrix} \in \mathbb{C}^{m_{tx} \times (m_{ue} + m_{tx})}
	\end{equation*}
	has a full Kruskal rank of $m_{tx}$. Assume $P$ is a Pareto-optimal precoder, such that all signals ${\cal S}_k(P) \neq 0$ are nonzero. Then either $P$ utilizes full power ((\ref{eq:powerrow}) are all fulfilled with equalities) or there exists $k$ such that
	\begin{equation}
	\frac{{\cal S}_k {\cal L}_k}{{\cal W}^2_k} \geq \frac{1}{m_{ue} - 1}. \label{eq:thcond}
	\end{equation}
	where ${\cal L}_k$ denotes 'leakage' according to a definition
	\begin{equation}
	{\cal L}_k := \sum \limits_{j \neq k}^{m_{ue}} |g_{jk}|^2. \label{eq:leakage}
	\end{equation}
\end{theorem}

Recall that a matrix has a full Kruskal rank if {\it each} (as compared to {\it one} for the full common rank) maximum-sized square submatrix of it is invertible. Assuming that the channel matrix values are commonly quite chaotic and never sparse (outside of Fourier domain), the used block channel-identity Kruskal rank assumption
\begin{equation*}
	\mbox{krank}(
	\begin{bmatrix}
		H & \vline & I_{m_{tx}}
	\end{bmatrix}
	) = m_{tx}
\end{equation*}
is quite realistic; particularly, it is satisfied with probability 1 if $H$ is random elementwise i.i.d. Gaussian.

The condition (\ref{eq:thcond}) serves as a formal description of the following concept: full power usage is {\it mandatory} for precoders that are Pareto-optimal for higher background noise coefficient values $\omega$, but is {\it not} a necessity in general. Further, in Section \ref{sec:exp}, we will provide empirical support for both of these claims, while showing that (\ref{eq:thcond}) is slightly pessimistic in practice.

This section of the paper will be devoted to proving Proposition \ref{th:power}. We will assume that for some precoder $P$ one of the per-antenna constraints (\ref{eq:powerrow}) is inactive; without loss of generality, let us assume it is the first one,
\begin{equation}
\label{eq:powfirst}
\| e_1^* P \|_2^2 < \beta_1.
\end{equation}
Then, we will search for a modification to the whole precoder
\begin{equation*}
P \rightarrow P + \epsilon \Delta P; \Delta P \in \mathbb{C}^{m_{tx} \times m_{ue}}, \epsilon > 0,
\end{equation*}
such that the power constraints for rows $2 \dots m_{tx}$ are not violated, and all the SINR values do not decrease, while at least one is strictly improved, and prove that (\ref{eq:thcond}) is then a necessity. In order to do so, we will utilize the gradients of the SINR values.

\subsection{SINR gradients}
Assuming $\omega_k > 0$, $\mathrm{SINR}_k$ is a smooth real-valued function of $m_{tx} \times m_{ue}$ complex parameters of $P$; in order to safely use differentiation, let us split the precoding matrix into a real and imaginary part, and denote
\begin{equation*}
	\frac{\partial \mathrm{SINR}_k}{\partial P} = \frac{\partial \mathrm{SINR}_k}{\partial \mbox{re}(P)} + i \frac{\partial \mathrm{SINR}_k}{\partial \mbox{im}(P)}
\end{equation*}
as the 'SINR gradient'. By direct real-valued differentiation, it follows that $\frac{\partial \mathrm{SINR}_k}{\partial P}$ has a closed-form representation as a rank-one matrixx:
\begin{align}
	\frac{\partial \mathrm{SINR}_k}{\partial P} & = h_k d_k^T, d_k \in \mathbb{C}^{m_{ue}} \label{eq:sinrderiv} \\
	(d_k)_j & = 
	\begin{cases}
		4 g_{kk} \frac{1}{{\cal W}_k}, j = k; \\	
		- 4 g_{kj} \frac{{\cal S}_k}{{\cal W}_k^2}, j \neq k.
	\end{cases} \nonumber
\end{align}
A rank-one structure is to be expected, since in the formula for one particular $\mathrm{SINR}_k$ each precoding vector is only used in a scalar product with $h_k$. As for $d_k$, the exact vector values have been carefully checked using numerical derivatives.

Utilizing the Frobenius norm induced scalar product
\begin{equation*}
\left( A, B \right)_F := tr(B^* A) = \sum \limits_{i,j} A_{ij} \overline{B_{ij}}; A, B \in \mathbb{C}^{m_{tx} \times m_{ue}}
\end{equation*}
and showing that a $\Delta P \in \mathbb{C}^{m_{tx} \times m_{ue}}$ exists that suffices the system of equations 
\begin{align*}
	(\Delta P, h_1 d_1^T)_F & = 1, &(\Delta P, e_2 e_2^* P)_F  & = -1, \\
	(\Delta P, h_2 d_2^T)_F & = 1, &(\Delta P, e_3 e_3^* P)_F  & = -1, \\
	& \dots & & \dots \\
	(\Delta P, h_{m_{ue}} d_{ue}^T)_F & = 1; &(\Delta P, e_{m_{tx}} e_{m_{tx}}^* P)_F &= -1;
\end{align*} 
would be {\it sufficient}\footnote{The equations enforce that the direction $\Delta P$ has a positive cosine with SINR derivatives and a negative cosine with power derivatives. They are sufficient but not necessary, of course: at the very least, the values of $\pm 1$ can be replaced with arbitrary positive and negative values respectively.} to show that the precoding matrix \mbox{$P + \epsilon \Delta P$} {\it strictly better} as compared to $P$ in the Pareto-sense for $\epsilon \rightarrow 0_+$. Vectorizing the unknown $\Delta P$, one could rewrite the system in a matrix form by replacing rank-one matrices with Kronecker products:
\begin{align}
	Y^* & \mbox{vec}(\Delta P) = \begin{bmatrix} 1 & \dots & 1 & -1 & \dots & -1 \end{bmatrix}^*, \label{eq:deltasys} \\
	Y & = \begin{bmatrix}
		h_1 & h_2 & \dots & h_{m_{ue}} & e_2 & e_3 & \dots & e_{m_{tx}} \\
		\otimes & \otimes & & \otimes & \otimes & \otimes & & \otimes \\
		d_1 & d_2 & \dots & d_{ue} & P^T e_2 & P^T e_3 & \dots & P^T e_{m_{tx}}
	\end{bmatrix}. \nonumber
\end{align}
Consequently, all that remains is to show that the long-column matrix \mbox{$Y \in \mathbb{C}^{m_{tx}m_{ue} \times (m_{tx} + m_{ue} - 1)}$} has a full column rank under the conditions of Proposition \ref{th:power}. Since it will be useful in the further analysis, let us introduce a matrix notation
\begin{equation*}
D^T := \begin{bmatrix}
	d_1 & d_2 & \dots & d_{ue}
\end{bmatrix} \in \mathbb{C}^{m_{ue} \times m_{ue}}
\end{equation*}
and establish a block formula for (\ref{eq:sinrderiv})
\begin{align}
D & = 4 \begin{bmatrix}
	\frac{1}{{\cal W}_1} & & \\
	& \ddots & \\
	& & \frac{1}{{\cal W}_{m_{ue}}}
\end{bmatrix}
G \circ F, \nonumber \\
F & :=
\begin{bmatrix}
	1 & - \frac{{\cal S}_1}{{\cal W}_1} & \dots & - \frac{{\cal S}_1}{{\cal W}_1} \\
	- \frac{{\cal S}_2}{{\cal W}_2} & 1 & \dots & - \frac{{\cal S}_2}{{\cal W}_2} \\
	\vdots & \vdots & \ddots & \vdots \\
	- \frac{{\cal S}_{m_{ue}}}{{\cal W}_{m_{ue}}} & - \frac{{\cal S}_{m_{ue}}}{{\cal W}_{m_{ue}}} & \dots & 1
\end{bmatrix}, \label{eq:deff}
\end{align}
where $G \circ F$ denotes Hadamard (elementwise) product.

\subsection{Rank of a Hadamard product}
The matrix $Y$ of interest has a peculiar structure - a set of columnwise Kronecker products. Bounding its rank {\it from below} is not so straightforward, but a work \cite{hadamard} provides an insightful tool for that - which requires a minor notation shift. Assume we denote two matrices of 'upper' and 'lower' Kronecker multiples, separated:
\begin{align*}
Y_{top} & \in \mathbb{C}^{m_{tx} \times (m_{tx} + m_{ue} - 1)}, Y_{bottom} \in \mathbb{C}^{m_{ue} \times (m_{tx} + m_{ue} - 1)}, \\
Y_{top} & := \begin{bmatrix}
	h_1 & \dots & h_{m_{ue}} & e_2 & e_3 & \dots & e_{m_{tx}}
\end{bmatrix}, \\
Y_{bottom} & := \begin{bmatrix}
	D^T & P^T e_2 & P^T e_3 & \dots & P^T e_{m_{tx}}
\end{bmatrix}.
\end{align*}
Then, we can switch to analyzing Hadamard product by using
\begin{equation*}
Y^* Y = Y_{top}^* Y_{top} \circ Y_{bottom}^* Y_{bottom}.
\end{equation*}
Both of the Hadamard multiples $Y_{top}^* Y_{top}$, $Y_{bottom}^* Y_{bottom}$ are nonnegative-definite, but are {\it not} full rank matrices; however, to show that $Y$ is full-rank, we still want $Y^* Y$ to be strictly {\it positive} definite. The work $\cite{hadamard}$ conveniently gives a tool exactly for that: it is sufficient to show that
\begin{equation*}
\mbox{krank}(Y_{top}) + \mbox{rank}(Y_{bottom}) > m_{tx} + m_{ue} - 1,
\end{equation*}
where the right side value is the number of columns of $Y, Y_{top}, Y_{bottom}$. Recalling the sizes, it can be seen that
\begin{equation*}
\mbox{krank}(Y_{top}) \leq m_{tx}, \mbox{rank}(Y_{bottom}) \leq m_{ue},
\end{equation*}
hence one of the matrices has to have a full Kruskal rank, and the other one should have full common rank. Since $Y_{top}$ is a submatrix of a block channel-identity matrix, it readily has a full Kruskal rank under the assumption of \mbox{Proposition \ref{th:power}}. Proving that $Y_{bottom}$ has a full common rank in all relevant cases is a much more complicated endeavor, however.

\subsection{Rank of an extended SINR derivative matrix}

Let us study the possible rank values of $Y_{bottom}$, considering four alternatives.

{\bf a)} The derivative submatrix is non-singular:
\begin{equation*}
\mbox{rank}(D) = m_{ue}.
\end{equation*}
This variant is what happens if a random $P$ is plugged in, and it guarantees $\mbox{rank}(Y_{bottom}) = m_{ue}$, which suits our purposes.

{\bf b)} The derivative submatrix is singular, but $Y_{bottom}$ has full rank, because some of the precoder rows number $2 \dots m_{tx}$ fill in the subspace complement. Again, since $\mbox{rank}(Y_{bottom}) = m_{ue}$, this variant suits our purposes.

{\bf c)} The derivative submatrix is singular, and the over-extended matrix 
\begin{equation*}
\begin{bmatrix}
D^T & P^T  
\end{bmatrix}
\xleftarrow{\text{up to perm.}} 
\begin{bmatrix}
Y_{bottom} & P^T e_1 
\end{bmatrix}
\end{equation*}
also doesn't have a full rank. This variant is actually impossible by construction unless one of the users is getting a zero signal; let's show it. The existence of a vector $d_{\perp} \neq 0$ such that $d_{\perp}^T Y_{bottom} = 0$ implies\footnote{due to lack of conjugation, $(\overline{d_{\perp}}, d_k)_2 = 0$, but $(d_{\perp}, d_k)_2 \neq 0$. The counter-intuitive notation is selected because otherwise the conjugation is present in every single formula afterwards and damages readability.} that
\begin{align}
0 = (D d_{\perp})_k & = ((G \circ F) d_{\perp})_k = g_{kk} d_{\perp, k} - \frac{{\cal S}_k}{{\cal W}_k} \sum \limits_{j \neq k} g_{kj} d_{\perp, j}; \nonumber \\ 
\frac{\sum \limits_{j \neq k} g_{kj} d_{\perp, j}}{g_{kk} d_{\perp, k}} & =  \frac{{\cal W}_k}{{\cal S}_k}. \mbox{   (if $g_{kk} d_{\perp, k} \neq 0$).} \label{eq:perpderiv}
\end{align}
However, using 
\begin{equation*}
P d_{\perp} = 0 \Rightarrow H^* P d_{\perp} = G d_{\perp} = 0,
\end{equation*}
we get
\begin{align*}
g_{kk} d_{\perp, k} & + \sum \limits_{j \neq k} g_{kj} d_{\perp, j} = 0, \\
\frac{\sum \limits_{j \neq k} g_{kj} d_{\perp, j}}{g_{kk} d_{\perp, k}} & = - 1, \mbox{ (if $g_{kk} d_{\perp, k} \neq 0$).}
\end{align*}
hence $g_{kk} d_{\perp, k} = 0, \forall k$, and since $d_{\perp} \neq 0$ - at least one of the users has a zero signal.

{\bf d)} The derivative submatrix is singular, $Y_{bottom}$ doesn't have a full rank, but the over-extended matrix
\begin{equation*}
	\begin{bmatrix}
		D^T & P^T  
	\end{bmatrix}
	\xleftarrow{\text{up to perm.}} 
	\begin{bmatrix}
		Y_{bottom} & P^T e_1 
	\end{bmatrix}
\end{equation*} 
has a full rank. That can actually happen in practice - such a combination $H, P, \omega$ will be presented later in Section \ref{sec:exp}. We will aim to show that this scenario is only possible for lower noises $\omega$. We could still follow the logic (\ref{eq:perpderiv}) of case {'c)'}, but the $P d_{\perp} = 0$ no longer holds - since $\overline{d_{\perp}}$ is no longer orthogonal to the first precoder row $P^T e_1$. 

\begin{lemma}{\it (Rank deficiency of $Y_{bottom}$):}
If all per-user signals ${\cal S}_k$ are nonzero, then
\begin{equation*}
\mbox{rank}(Y_{bottom}) \geq m_{ue} - 1.
\end{equation*}
\end{lemma}
\begin{IEEEproof}
Assume again that $\overline{d_{\perp}}$ is orthogonal to $d_1 \dots d_{ue}$, then (\ref{eq:perpderiv}) holds. We can see that $P d_{\perp} = \gamma e_1, \gamma \in \mathbb{C} \setminus \{0\}$, since $Y_{bottom}$ contains all the transposed rows of $P$ except the first one. Hence,
\begin{align*}
G d_{\perp} = H^* P d_{\perp} & = \gamma H^* e_1, \\
g_{kk} d_{\perp, k} + \sum \limits_{j \neq k} g_{kj} d_{\perp, j} & = \gamma h_k^* e_1, \\
\frac{{\cal S}_k}{{\cal W}_k} g_{kk} d_{\perp, k} + \frac{{\cal S}_k}{{\cal W}_k} \sum \limits_{j \neq k} g_{kj} d_{\perp, j} & = \gamma h_k^* e_1, \\
\mbox{(summing with (\ref{eq:perpderiv})) } (\frac{{\cal S}_k}{{\cal W}_k} + 1) g_{kk} d_{\perp, k} & = \gamma \frac{{\cal S}_k}{{\cal W}_k}  h_k^* e_1,
\end{align*}
and hence all the possible $d_{\perp}$ such that $d_{\perp}^T Y_{bottom}$ are collinear with a one-dimensional direction
\begin{align}
& d_{\perp} \parallel H^* e_1 \circ v; \label{eq:dperpcollin} \\
v & := \begin{bmatrix} \frac{1}{(1 + \frac{{\cal W}_1}{{\cal S}_1}) g_{11}} & \frac{1}{(1 + \frac{{\cal W}_2}{{\cal S}_2}) g_{22}} & \dots & \end{bmatrix}^T. \nonumber
\end{align}
\end{IEEEproof}

The next Lemma will utilize the explicit direction formula (\ref{eq:dperpcollin}) to construct a concise necessary condition for a rank-deficient $Y_{bottom}$.

\begin{lemma}{\it (Unit eigenvalue):}
\label{lemma:unit}
	If all per-user signals ${\cal S}_k$ are nonzero, and $Y_{bottom}$ is rank-deficient, then the matrix
	\begin{equation*}
	G \mathrm{diag}(v) \mathrm{\mbox{   } (see \mbox{ } (\ref{eq:dperpcollin}))}
	\end{equation*}
	has an eigenvalue equal to 1.
\end{lemma}
\begin{IEEEproof}
Plugging the direction of $d_{\perp}$, established in the previous Lemma, into $D d_{\perp} = 0$, and rewriting the vector Hadamard product using a diagonal matrix gives
\begin{equation*}
(G \circ F) (H^* e_1 \circ v) = (G \circ F) \mathrm{diag}(v) H^* e_1 = 0.
\end{equation*}
Let us construct a way to open a matrix Hadamard product utilizing the structure of matrix F. Recalling (\ref{eq:deff}) one can utilize
\begin{align} 
F & = F_d - F_1; \nonumber \\
F_d & := \mathrm{diag}(1 + \frac{{\cal S}_1}{{\cal W}_1}, 1 + \frac{{\cal S}_2}{{\cal W}_2} \mbox{ } \dots \mbox{ } 1 + \frac{{\cal S}_{m_{ue}}}{{\cal W}_{m_{ue}}}) \nonumber \\
F_1 & :=  \begin{bmatrix}
	\frac{{\cal S}_1}{{\cal W}_1} \\
	\frac{{\cal S}_2}{{\cal W}_2} \\
	\vdots \\
	\frac{{\cal S}_{m_{ue}}}{{\cal W}_{m_{ue}}}
\end{bmatrix}
\begin{bmatrix}
	1 & 1 & \dots & 1 
\end{bmatrix}, \nonumber \\
(G & \circ F_1) \mathrm{diag}(v) H^* e_1 = (G \circ F_d) \mathrm{diag}(v) H^* e_1. \label{eq:fsplit}
\end{align}
Both Hadamard products in the last equation can be simplified, because $F_d$ is diagonal, and $F_1$ is rank-one. Indeed, $(G \circ F_d) \mathrm{diag}(v)$ is diagonal: direct multiplication gives that
\begin{equation*}
((G \circ F_d)\mathrm{diag}(v))_{kk} = \frac{g_{kk} (1 + \frac{{\cal S}_k}{{\cal W}_k})}{(1 + \frac{{\cal W}_k}{{\cal S}_k}) g_{kk}} = \frac{{\cal S}_k}{{\cal W}_k}.
\end{equation*}
A Hamamard product with a rank-one matrix corresponds to a left and right common multiplication by a diagonalization of rank-one factors, hence
\begin{equation*}
(G \circ F_1) \mathrm{diag}(v) = \mathrm{diag}(\frac{{\cal S}_1}{{\cal W}_1}, \frac{{\cal S}_2}{{\cal W}_2} \mbox{ }\dots\mbox{ }) G \mathrm{diag}(v).
\end{equation*}
Revisiting (\ref{eq:fsplit}) then gives
\begin{equation*}
G \mathrm{diag}(v) H^* e_1 = H^* e_1.
\end{equation*}
\end{IEEEproof}

Lemma \ref{lemma:unit} finally gives a clear path towards the proof of Proposition \ref{th:power}. To see that, let us draw the product $G \mathrm{diag}(v)$:
\begin{equation*}
\begin{bmatrix}
\frac{1}{1 + \frac{{\cal W}_1}{{\cal S}_1}} \mbox{  } & \frac{1}{1 + \frac{{\cal W}_2}{{\cal S}_2}} \frac{g_{12}}{g_{22}} & \dots &  \frac{1}{1 + \frac{{\cal W}_{m_{ue}}}{{\cal S}_{m_{ue}}}} \frac{g_{1 m_{ue}}}{g_{m_{ue} m_{ue}}} \\
\frac{1}{1 + \frac{{\cal W}_1}{{\cal S}_1}} \frac{g_{21}}{g_{11}} & \frac{1}{1 + \frac{{\cal W}_2}{{\cal S}_2}} \mbox{  } & \dots &  \frac{1}{1 + \frac{{\cal W}_{m_{ue}}}{{\cal S}_{m_{ue}}}} \frac{g_{2 m_{ue}}}{g_{m_{ue} m_{ue}}} \\
\vdots & \vdots & \ddots & \vdots \\
\frac{1}{1 + \frac{{\cal W}_1}{{\cal S}_1}} \frac{g_{m_{ue} 1}}{g_{11}} & \frac{1}{1 + \frac{{\cal W}_2}{{\cal S}_2}} \frac{g_{m_{ue} 2}}{g_{22}} & \dots & \frac{1}{1 + \frac{{\cal W}_{m_{ue}}}{{\cal S}_{m_{ue}}}} \mbox{   }
\end{bmatrix}
\end{equation*}
Note that the diagonal elements are real values bounded by one, which only approach the limit if the SINR values $\frac{{\cal S}_k}{{\cal W}_k}$ are large, and the off-diagonal elements correspond to scaled relative interference/leakage. Proposition \ref{th:power} is then constructed by assessing when the matrix above could have an eigenvalue equal to one by utilizing Gershgorin circles: such an assessment results in a more concise condition for columns, but then utilizes the concept of leakage ${\cal L}_k$ defined in (\ref{eq:leakage}).

\begin{IEEEproof}[Proof of Proposition \ref{th:power}]
Consider column-wise Gershgorin circles as a bound for eigenvalues of $G \mathrm{diag}(V)$. Following Lemma \ref{lemma:unit}, if $P$ doesn't utilize full power, it is then necessary that for some column $k$
\begin{align*}
	\frac{1}{1 + \frac{{\cal W}_k}{{\cal S}_k}} (1 & + \sum \limits_{j \neq k}  \frac{\lvert g_{jk} \rvert }{\lvert g_{kk} \rvert}  ) \geq 1; \\
	1 + \sum \limits_{j \neq k} \frac{\lvert g_{jk} \rvert}{ \lvert g_{kk} \rvert } & \geq 1 + \frac{{\cal W}_k}{{\cal S}_k}; \\
	(m_{ue} - 1) \frac{{\cal L}_k}{{\cal S}_k} & \geq (\sum \limits_{j \neq k} \frac{\lvert g_{jk} \rvert}{ \lvert g_{kk} \rvert })^2 \geq \frac{{\cal W}^2_k}{{\cal S}^2_k}.
\end{align*}
\end{IEEEproof}

\begin{remark}{{\it (Zero signals)}}
The condition of all signals ${\cal S}_k$ being nonzero is not crucial for Proposition \ref{th:power}. Indeed, if ${\cal S}_1(P) = 0$, it is sufficient to see that
\begin{equation*}
	\widetilde P = \begin{bmatrix}
		0 & P_{shrunk} \\
	\end{bmatrix}, P_{shrunk} := \begin{bmatrix}
	p_2 & p_3 & \dots & p_{m_{ue}} \\
	\end{bmatrix}
\end{equation*}
cannot utilize more power than the initial $P$ and then one can apply Proposition \ref{th:power} to $P_{shrunk}$ and
\begin{align*}
	H_{shrunk} & := \begin{bmatrix}
		h_2 & h_3 & \dots & h_{m_{ue}} \\
	\end{bmatrix}, \\
	H_{shrunk},& P_{shrunk}  \in \mathbb{C}^{m_{tx} \times (m_{ue} - 1)}.
\end{align*}
One would then have to plug ${\cal S}(P_{shrunk})$, ${\cal W}(P_{shrunk})$, ${\cal L}(P_{shrunk})$ into the condition (\ref{eq:thcond}), however. 
\end{remark}

\section{Algorithm for building Pareto-optimal precoders}
\subsection{Iterative refinement}
\label{subsec:refine}
Let us now discuss how the full power necessity helps building practical precoding matrices. Recall the introduction on the Quality-of-Service problem in Subsection \ref{subsec:qos}: using convex optimization framework, the previous work \cite{bjornsson} establishes that every single Pareto-optimal SINR combination can be obtained using a parameterization (\ref{eq:noncomp}). This parameterization is effectively based on $m_{tx} + m_{ue} - 1$ real Lagrange multiplier parameters, but this system of parameters is redundant, since the Pareto boundary itself should be a surface in $\mathbb{R}^{m_{ue}}$ with at most $m_{ue} - 1$ parameters.

Let us return to the idea of full power necessity - and try to impose $m_{tx}$ 'extra' full per-antenna power constraints. Assume $\lambda_k$ are the 'true' parameterization parameters, and\footnote{By (\ref{eq:noncomp}) the directions $\hat p_j$ are self-scaled, hence the parameterization is invariant under uniform scaling of $\lambda, \mu$.}
\begin{equation*}
\sum \limits_{k = 1}^{m_{ue}} \lambda_k = 1.
\end{equation*} 
Let us attempt to find $\mu = \mu(\lambda)$ using the following system of equations:
\begin{equation*}
\|e^*_i P [\lambda, \mu(\lambda)] \|_2^2 = \beta_i,\mbox{ } i \in \{1, \dots, m_{tx}\}.
\end{equation*} 

Empirically, it turns out a simple iterative numerical scheme, that rescales each $\mu_i(\lambda)$ based on how close $\|e^*_i P [\lambda, \mu(\lambda)] \|_2^2$ was to $\beta_i$ in the previous iteration, converges very well:

\begin{itemize}
	\item Set an initial guess $\mu^{(0)}_{i} = \frac{1}{m_{tx}}$. Following the convex QoS derivation a precoder $P(\lambda, \mu^{(0)})$ is actually {\it Pareto-optimal} for Frobenius power constraints (\ref{eq:powerfro}) \cite{bjornsson}.
	\item Set iteration counter $n = 0$ and repeat
		\begin{itemize}
		\item Compute $P(\lambda, \mu^{(n)})$ and the per-row power consumption rates
		\begin{equation*}
		\alpha_i^{(n)} := \frac{\| e_i^* P(\lambda, \mu^{(n)}) \|_2}{\sqrt{\beta_i}} 
		\end{equation*}
		\item Update and rescale 
		\begin{equation}
		\mu^{(n + 1)}_i := \frac{\mu^{(n)}_i \alpha_i}{\sum \limits_{k = 1}^{m_{tx}} \mu^{(n)}_i \alpha_i}, \label{eq:muupdate}
		\end{equation} then increase iteration counter $n := n + 1$.
		\end{itemize}
	until, for a tolerance parameter $\delta$,
	\begin{equation}
	\min  \limits_{i = 1 \dots m_{tx}} \alpha^{(n)}_{i} > 1 - \delta, \max \limits_{i = 1 \dots m_{tx}} \alpha^{(n)}_{i} < \frac{1}{1 - \delta}. \label{eq:muexit}
	\end{equation}
	\item Rescale 
	\begin{equation}
	P(\lambda, \mu_{\delta}(\lambda)) := (1 - \delta) P(\lambda, \mu^{(n)}). \label{eq:muscale}
	\end{equation}
\end{itemize}

Example values of $\delta$ used in experiments were in range $[10^{-2}, 10^{-8}]$. Following the cycle exit condition, the precoding matrix 'almost' satisfies the per-antenna power constraints with a possible excess defined by $\delta$, and the minor scaling (\ref{eq:muscale}) guarantees that the final precoder is entirely fair.

It turns out the proposed numerical scheme has a very stable convergence empirically - moreover, the convergence speed in terms of sufficient iterations $n$ does not really depend on the {\it sizes} of the precoding matrix, but only on the tolerance $\delta$. Table \ref{tab:iter} provides suffcient-to-convergence times $P(\lambda, \mu)$ is computed according to formula (\ref{eq:noncomp}) - these values were averaged over $10000$ experiments with random $\lambda$, each value of which was selected randomly using uniform distribution on $[0, 1]$ and post-scaled towards a sum of $1$.

%eAdditionally, an extreme case was considered, in which only one random value of $\lambda_k$ was set to $1$, and the remaining values were set to zero; it turned out that in this case the number of iterations might grow by a factor of $2$; however this might not be a relevant value since if most user signal amplitudes are preset to zero from the beginning, then the precoder may be effectively computed in a space $P_{shrunk} \in \mathbb{C}^{m_{tx} \times 1}$, compensating the added iterations in complexity.

From Table \ref{tab:iter}, it is clear that the iteration number growth does not really depend on the sizes of $H, P$, and depends logarithmically on $\frac{1}{\delta}$. In the next subsection we will argue that the formulas (\ref{eq:noncomp}) should not be used directly - because there is a way to rewrite those with far better complexity.

%In the numerical experiment section we will show that running the algorithm till convergence is mostly interesting for theoretical 'upper-bound' studies, while for practical application basically one or two computations of $P(\lambda, \mu)$ is sufficient to obtain highly competitive solutions. 

\begin{table}
\centering
\caption{Average sufficient iteration numbers and corresponding standard deviations in dependence on $\delta, m_{tx}, m_{ue}$.}
\label{tab:iter}
\begin{tabular}{|c|c|c|c|}
\hline
\rule{0pt}{14pt} $\delta$ & $H, P \in \mathbb{C}^{8 \times 2}$ & $H, P \in \mathbb{C}^{24 \times 8}$ & $H, P \in \mathbb{C}^{192 \times 24}$ \\
\hline
\rule{0pt}{9pt} $10^{-2}$ & $5.21 \pm 1.38$ & $4.58 \pm 0.50$ & $5.15 \pm 0.41$ \\
\hline
\rule{0pt}{9pt} $10^{-4}$ & $12.02 \pm 1.58$ & $10.10 \pm 0.79$ & $10.48 \pm 0.56$ \\
\hline
\rule{0pt}{9pt} $10^{-6}$ & $19.04 \pm 1.78$ & $15.99 \pm 1.21$ & $15.91 \pm 0.66$ \\
\hline
\rule{0pt}{9pt} $10^{-8}$ & $26.22 \pm 2.02$ & $22.02 \pm 1.61$ & $21.39 \pm 0.77$ \\
\hline
\end{tabular}
\end{table}

\subsection{Efficient parameterization implementation}
Let us now focus on a numerical simplification of (\ref{eq:noncomp}). We are going to utilize a Woodbury matrix identity \cite{woodbury}
\begin{align}
	(A + B B^*)^{-1} & = A^{-1} - A^{-1} B (I + B^* A^{-1} B)^{-1} B^* A,  \nonumber \\
	A \in \mathbb{C}^{m \times m}, B & \in \mathbb{C}^{m \times r}, \det(I + B^* A^{-1} B) \neq 0. \label{eq:rankupd}
\end{align}
to process low-rank updates of an inverse multiple times. Let us denote
\begin{equation*}
	\widetilde H =
	H \mathrm{diag}(\frac{\sqrt{\lambda_1}}{\omega_1} \dots \frac{\sqrt{\lambda_{m_{ue}}}}{\omega_{m_{ue}}}).
	%\begin{bmatrix}
	%\frac{\sqrt{\lambda_1}}{\omega_1} & & \\
	%& \ddots & \\
	%& & \frac{\sqrt{\lambda_{m_{ue}}}}{\omega_{m_{ue}}}
	%\end{bmatrix}
\end{equation*}
Then, recalling (\ref{eq:noncomp}), note that the values of $\hat p_j$ and the values of signal and interference coefficients $S$ are defined by two matrices
\begin{align}
(M + \widetilde H \widetilde H^*)^{-1} \widetilde H & \in \mathbb{C}^{m_{tx} \times m_{ue}} \label{eq:cow} \\
Z := \widetilde H^* (M + \widetilde H \widetilde H^*)^{-1} & \in \mathbb{C}^{m_{ue} \times m_{ue}} \nonumber \\
\mbox{where } M := \mathrm{diag}(\frac{\mu_1}{\beta_1} \dots \frac{\mu_{m_{tx}}}{\beta_{m_{tx}}}) & \in \mathbb{C}^{m_{tx} \times m_{tx}}.
%M & := \begin{bmatrix}
%\frac{\mu_1}{\beta_1} & & \\
%& \ddots & \\
%& & \frac{\mu_{m_{tx}}}{\beta_{m_{tx}}}
%\end{bmatrix}. \nonumber
\end{align}
We can utilize (\ref{eq:rankupd}) to open $(M + \widetilde H \widetilde H^*)^{-1}$ with $r = \mbox{rank}(\widetilde H \widetilde H^*) = m_{ue}$ and establish\footnote{If $C > 0$, then using diagonalization, one can see $C - C (I + C)^{-1} C > 0$.}
\begin{align}
	\widetilde H^* (M + \widetilde H \widetilde H^*)^{-1} \widetilde H & = C  - C (I + C)^{-1} C \label{eq:midprec}  \\
	(M + \widetilde H \widetilde H^*)^{-1} \widetilde H & = M^{-1} \widetilde H (I - (I + C)^{-1} C), \nonumber
\end{align}
where $C := \widetilde H^* M^{-1} \widetilde H$. Furthermore, using (\ref{eq:rankupd}) with $r = 1$ and appropriate signs also gives 
\begin{align}
	\gamma_j  &= \tilde h_j^* (\Psi - \tilde h_j \tilde h_j^*)^{-1} \tilde h_j \nonumber \\
	& = \tilde h_j^* \Psi^{-1} \tilde h_j + \tilde h_j^* \Psi^{-1} \tilde h_j (1 - \tilde h_j^* \Psi^{-1} \tilde h_j)^{-1} \tilde h_j^* \Psi^{-1} \tilde h_j \nonumber \\
	& = z_{jj} + \frac{z_{jj^2}}{1 - z_{jj}} = \frac{z_{jj}}{1 - z_{jj}}, \label{eq:gamma}
\end{align}
where the last transformation follows from $h_j^* \Psi^{-1} \tilde h_j$ being a diagonal element of $Z$ defined in (\ref{eq:cow}).
Plugging (\ref{eq:midprec}) and (\ref{eq:gamma}) into (\ref{eq:noncomp}) results in Algorithm \ref{alg:twopar}, which has a complexity of $O(m_{tx} m_{ue}^2)$. The two dominant operations are labelled with $\clubsuit$ and have the same complexity constant as the dominant operations of the common Zero-Forcing implementations: computing the covariance and multiplying the channel by a square matrix - but in Algorithm \ref{alg:twopar} the channel is preprocessed by a Hadamard product with a rank-one matrix.  
The user-cubic $O(m_{ue}^3)$ part of the complexity brings a larger constant as compared to Zero-Forcing; however, preemptive tridiagonalization of positive-definite $C$ can be utilized to effectively simplify the operations marked with $\spadesuit$.

\begin{algorithm}[h!]
	\caption{\it Parametric precoding $P(\lambda, \mu)$}
	\label{alg:twopar}
	\KwInput{Channel $H$, power constraints $\beta_j > 0$, noises $\omega_k > 0$, weights $\lambda_j > 0, \sum \limits_{j = 1}^{m_{ue}} \lambda_j = 1.$ and $\mu_i > 0, \sum \limits_{i = 1}^{m_{tx}} \mu_i = 1.$}
	$\widehat H := \mathrm{diag}(\sqrt{\frac{\beta_1}{\mu_1}} \dots \sqrt{\frac{\beta_{m_{tx}}}{\mu_{m_{tx}}}}) \mbox{ $H$ } \mathrm{diag}(\frac{\sqrt{\lambda_1}}{\omega_1} \dots \frac{\sqrt{\lambda_{m_{ue}}}}{\omega_{m_{ue}}})$ \;
	$C := \widehat H^* \widehat H$ \mbox{ $\clubsuit$} \;
	$\widehat P := \mathrm{diag}(\sqrt{\frac{\beta_1}{\mu_1}} \dots \sqrt{\frac{\beta_{m_{tx}}}{\mu_{m_{tx}}}}) \mbox{ $\widehat H$ } (I - (I + C)^{-1} C)$ \mbox{ $\clubsuit \spadesuit$} \;
	$Z := C - C (I + C)^{-1} C \mbox{ $\spadesuit$}; \gamma_j := \frac{Z_{jj}}{1 - Z_{jj}}, j = 1 \dots m_{ue}$ \;
	$S_{kj} :=
	\begin{cases}
		\lvert Z_{kj} \rvert^2 / \| p_k \|_2^2, i = j \\
		- \gamma_j \lvert Z_{kj} \rvert^2 / \| p_k \|_2^2, k \neq j.
	\end{cases}, S \in \mathbb{R}^{m_{ue} \times m_{ue}}$\;
	$\kappa := S^{-T} 
	\begin{bmatrix}
		\gamma_1 \lambda_1 & \dots & \gamma_{m_{ue}} \lambda_{m_{ue}}
	\end{bmatrix}^T; p_j := \hat p_j \frac{\sqrt{\kappa_j}}{\|\vec p_j \|_2}.$
\end{algorithm}

\begin{remark}{{\it (Numerical instability)}}
\label{rem:instability}
Assuming $C$ is diagonal and $\omega \rightarrow 0_+$ in the Algorithm \ref{alg:twopar} brings us to two divisions
\begin{equation*}
f_j := \frac{c_j}{1 + c_j}; \gamma_j := \frac{f_j}{1 - f_j}.
\end{equation*}
Since $C$ is rescaled by the inverses of $\omega$, $c_j$ is a potentially large real value - and the above two operations are dangerous in double precision with SINR values above $10^8$. To prevent $\gamma_j$ showing up as negative because of the rounding errors in $f_j \approx 1$, either avoid extreme scenarios or use higher precision.
\end{remark}

\begin{figure*}
	\centering
	\subfloat[Parametric precoder with equal $\mu_j$]{\includegraphics[width=0.3\linewidth]{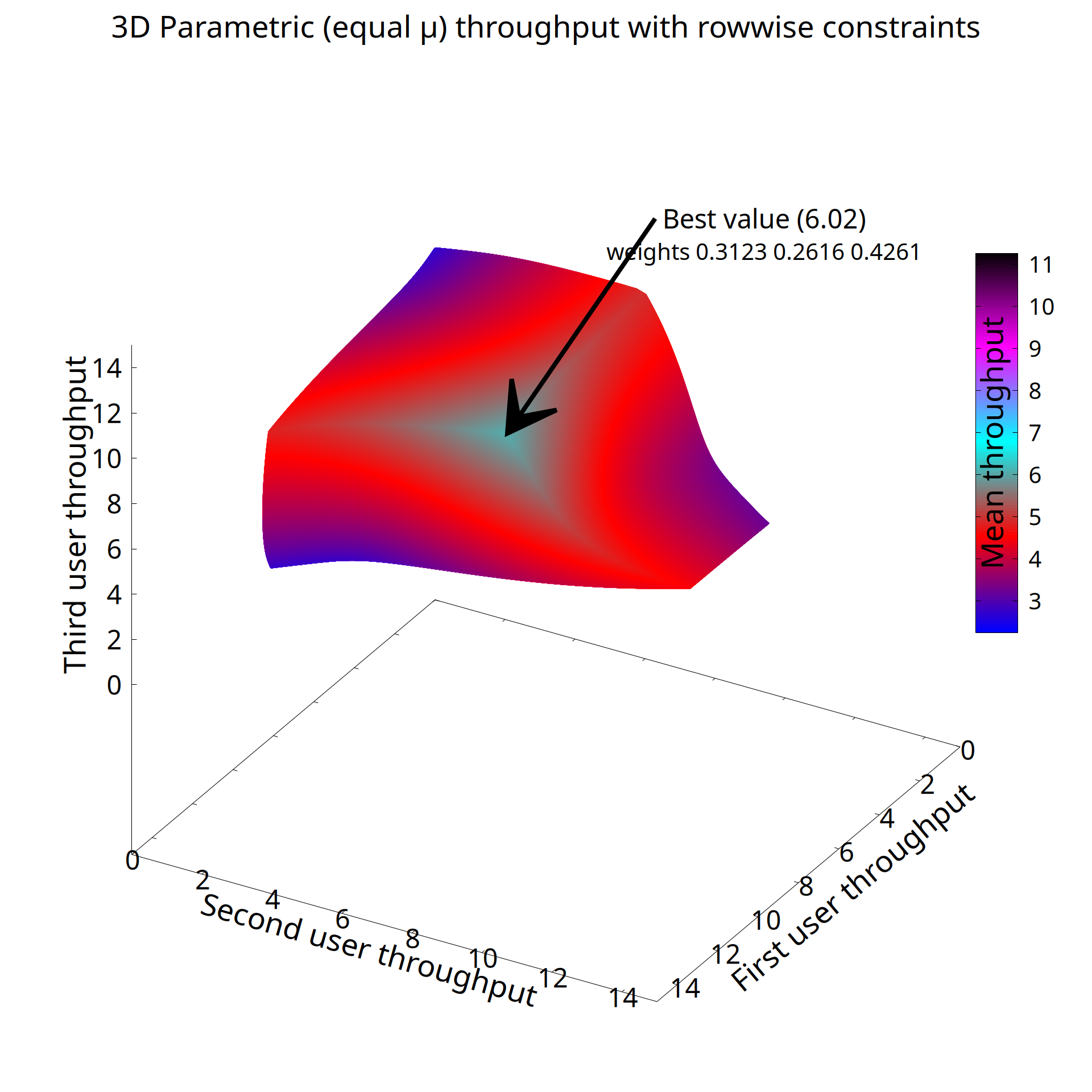}
		\label{subfig:toy_p0}}
	\hfil
	\subfloat[Parametric precoder, one iteration over $\mu$]{\includegraphics[width=0.3\linewidth]{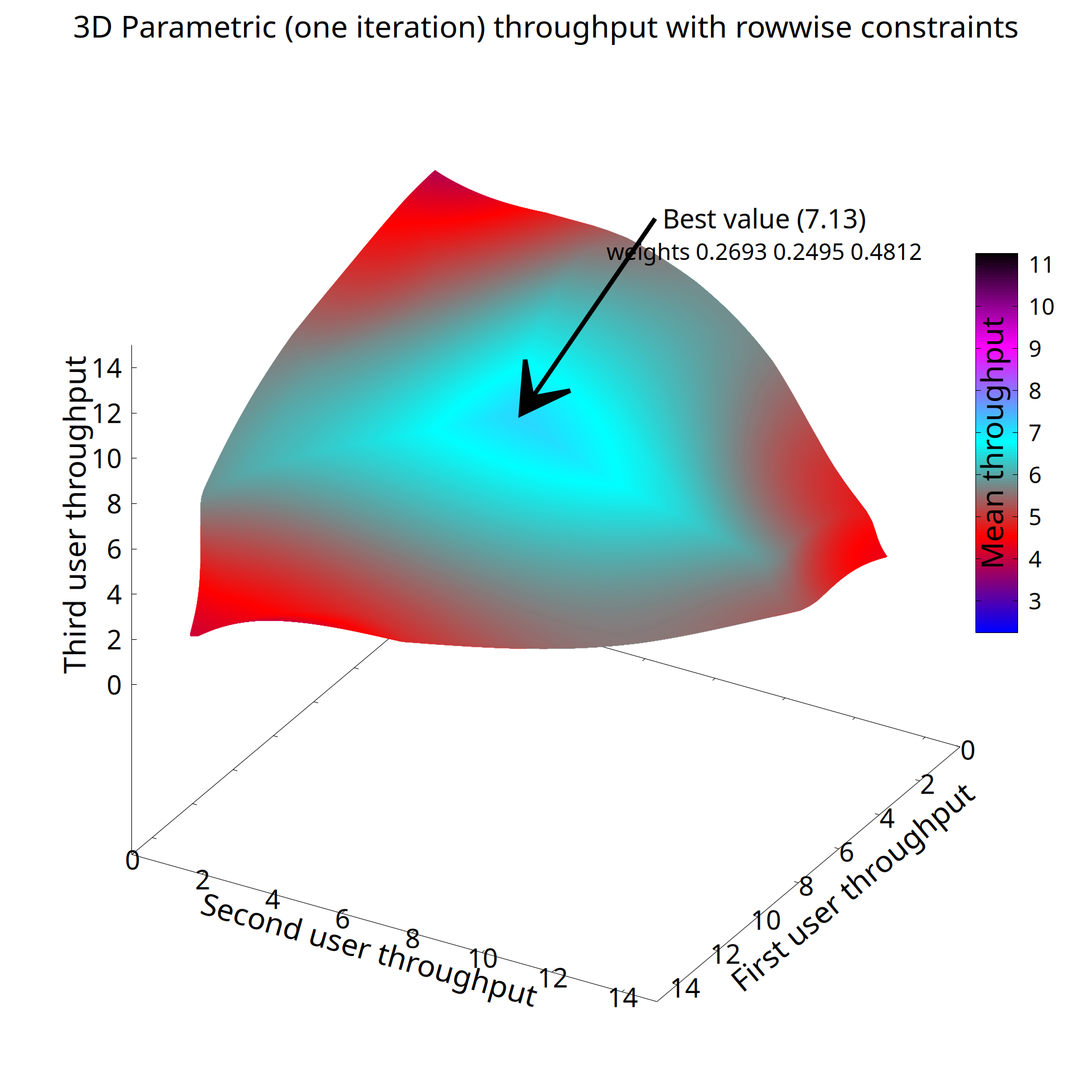}
		\label{subfig:toy_p1}}
	\hfil
	\subfloat[Parametric precoder, convergent to $\delta = 0.01$ ]{\includegraphics[width=0.3\linewidth]{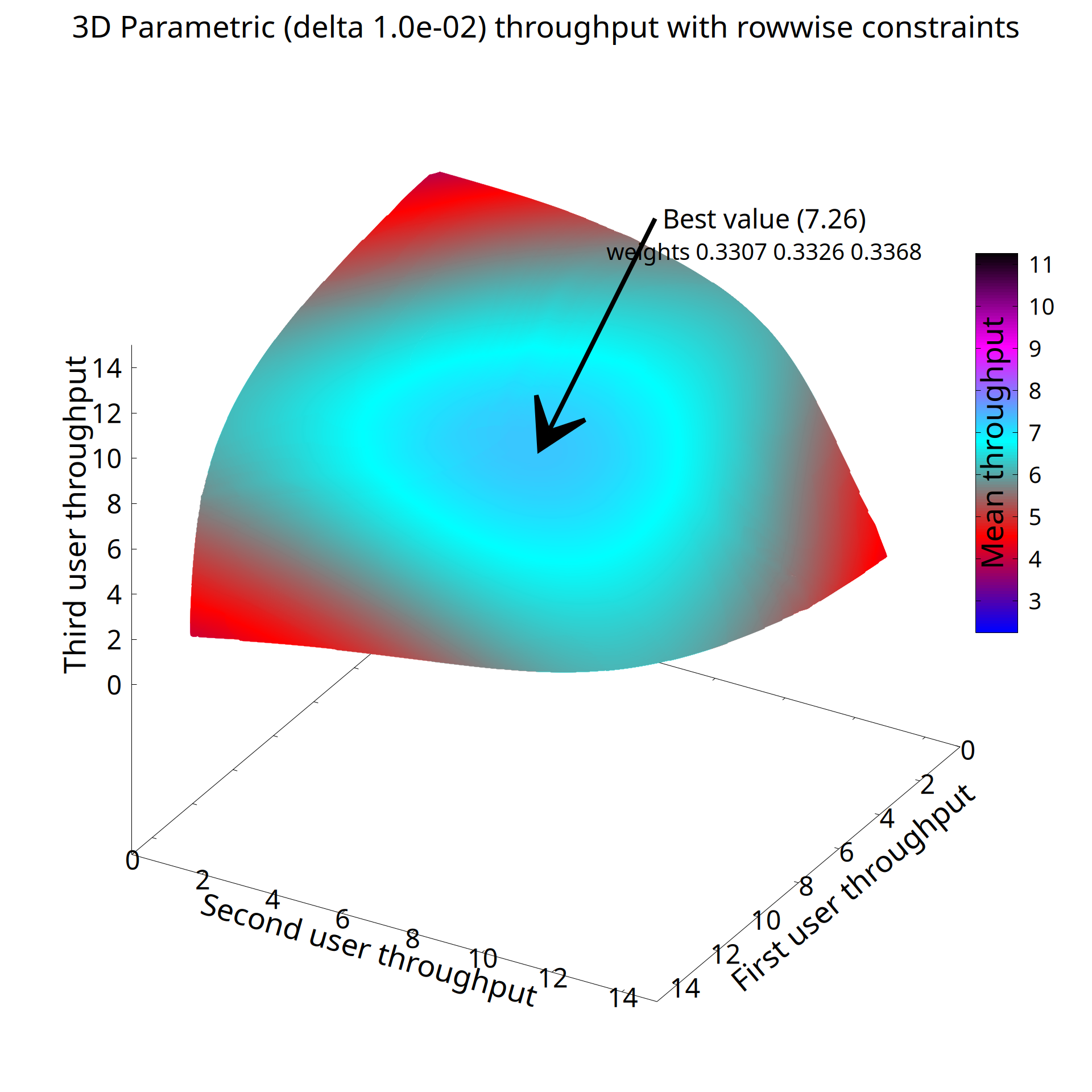}
		\label{subfig:toy_p5}}
		
	\subfloat[Channel values, problem parameters and ZF/SLNR references]
	{
		\begin{minipage}{0.32\linewidth}
		$H = \begin{bmatrix}
			-0.4 & -0.2 & 0.0 \\
			-0.4 & 0.9 & 1.2 \\
			-0.2 & 0.1 & -0.2 \\
			0.7 & -0.8 & -0.5 \\
			0.0 & -0.2 & -0.2 \\
			0.6 & -0.5 & 0.6 \\
			-0.4 & 1.2 & -0.1 \\
			-1.2 & -0.4 & -0.8 
		\end{bmatrix} \in \mathbb{R}^{8 \times 3}.$
		\begin{center}
		$\omega_k = 1, \beta_i = 1.$
		\end{center}
		\begin{center}
		\centering $k \in \{1, \dots, m_{ue}\}; \mbox{ } i \in \{1, \dots, m_{tx}\}$
		\end{center}
		\end{minipage}
		\hfil
		\begin{minipage}{0.32\linewidth}
		\includegraphics[width=1.0\linewidth]{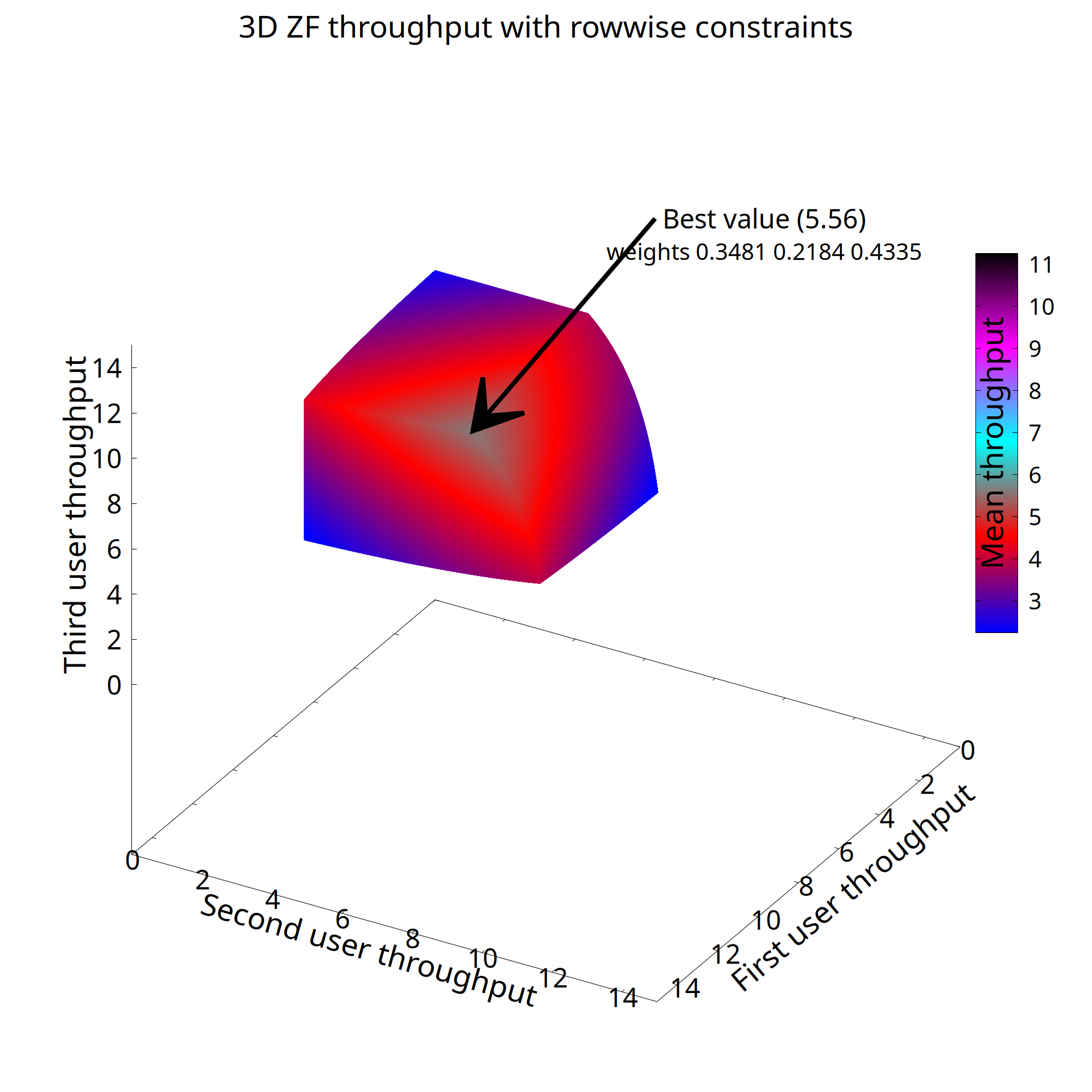}
		\end{minipage}
		\hfil
		\begin{minipage}{0.32\linewidth}
		\includegraphics[width=1.0\linewidth]{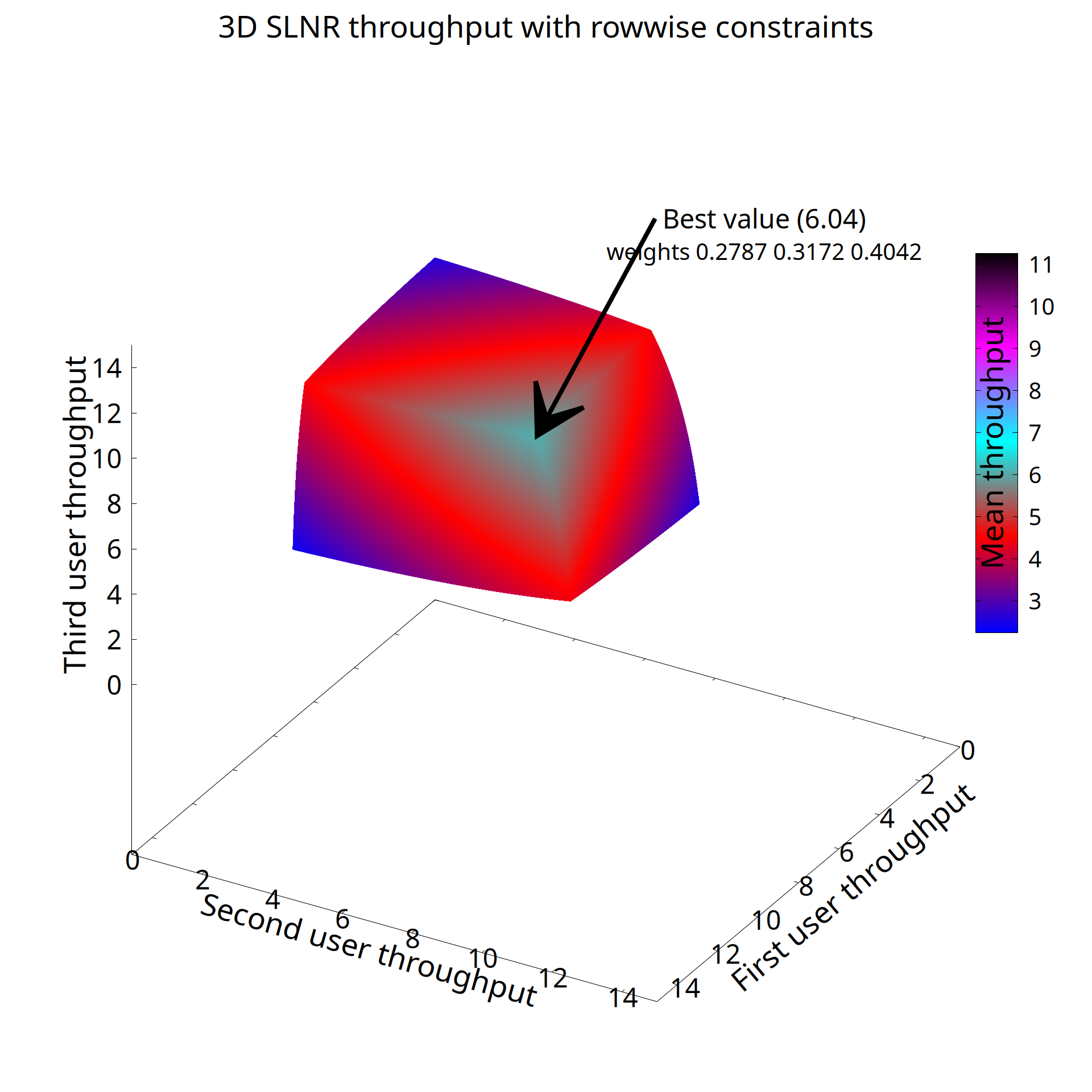}
		\end{minipage}
		\label{subfig:toy_ch_zf_slnr}
	}
	
	%\centering
	%\subfloat[Zero-forcing reference]{\includegraphics[width=0.3\linewidth]{pics/toy/snapshot3D_p2.png}
	%	\label{subfig:toy_zf}}
	%\hfil
	%\subfloat[SLNR reference]{
	%	\includegraphics[width=0.3\linewidth]{pics/toy/snapshot3D_p1.png}
	%	\label{subfig:toy_slnr}}
	
	\caption{Simulation results on a toy channel of size $H \in \mathbb{R}^{8 \times 3}$. Each axis represents a value $10 \log_{10} (1 + \mathrm{SINR}_k), k = 1,2,3.$ Color shows the mean value. The graphs are constructed pointwise using $1.25 \times 10^7$ points, and the arrows show the point with the highest mean throughput.}
	\label{fig:toy}
\end{figure*}

\section{Numerical experiments}
\label{sec:exp}

\subsection{Reference methods}
Let us provide a numerical study of the proposed algorithm framework of combining Algorithm \ref{alg:twopar} combined with an iterative refinement of Subsection \ref{subsec:refine}. Since we were aiming at studying the Pareto surfaces of the SINR multi-objective. Let us describe our references. We will use an SVD-notation
\begin{equation}
\label{eq:svd}
H = U \Sigma V^*, U \in \mathbb{C}^{m_{tx} \times m_{ue}}, V \in \mathbb{C}^{m_{ue} \times m_{ue}},
\end{equation}
for reference description purposes.

\begin{itemize}
\item {\it Zero-Forcing reference.} Arguably the most well-known precoding algorithm reference, a Zero-Forcing solution we denoted previously as $P_{\mathrm{ZF}} = H (H^* H)^{-1}$ in the introduction is now easier to perceive as
\begin{equation*}
P_{\mathrm{ZF}} = U \Sigma^{-1} V^*.
\end{equation*}
This solution manifests a zero-interference idea \mbox{$H^* P_{\mathrm{ZF}} = I$}, and is commonly post-processed by
\begin{align}
P_{\mathrm{ZF}, pwr}(\kappa) & = 
\frac{P_{\mathrm{ZF}} D_{\kappa}} 
{\max \limits_{i} \left( \frac{1}{\sqrt{\beta_i}}\| e_i^* P_{\mathrm{ZF}}
	D_{\kappa} \|_2 \right)
}, \label{eq:pwr1} \\
%D_{\kappa} & := \begin{bmatrix}
%	\sqrt{\kappa_1} & & \\
%	& \ddots & \\
%	& & \sqrt{\kappa_{m_{ue}}}
%\end{bmatrix}
D_{\kappa} & = \mathrm{diag}(\sqrt{\kappa_1}, \mbox{ } \dots \mbox{ }, \sqrt{\kappa_{m_{ue}}}) \nonumber
\end{align}
with various weights\footnote{Surely, we will not be able to check {\it all} possible 'power allocations' of $\kappa$ since they form a continuum, but we can only build readable surface graphs for $m_{ue} \leq 3$; for this dimensionality a large enough discrete set of vectors $\kappa \in \mathbb{R}_+^{m_{ue}}$ will suffice for demonstrations.} $\kappa \geq 0$ limited by \mbox{$\sum \limits_{j = 1}^{m_{ue}} \kappa_j = 1$}.

\item {\it SLNR maximizer reference.} Another known precoding algorithm in the form
\begin{equation}
	P_{\mathrm{SLNR}} = U W^*, \mbox{ } W_{kj} := \frac{\sigma_j}{\omega_k^2 + \sigma_j^2} V_{kj}, \label{eq:slnr}
\end{equation}
($\sigma_j$ are the singular values from $\Sigma$) is constructed by optimizing a surrogate for SINR:
\begin{equation*}
\mathrm{SLNR}_k := \frac{{\cal S}_k}{{\cal L}_k + {\cal N}_k} = \frac{|h_k^* p_k|^2}{\omega_k^2 + \sum \limits_{j = 1, j \neq k}^{m_{ue}} |h_j^* p_k|^2},
\end{equation*}
where SLNR stands for 'Signal-to-Leakage-and-Noise-Ratio'  \cite{slnr1, slnr2}. Again, we will utilize 'power allocated'
\begin{equation}
	P_{\mathrm{SLNR}, pwr}(\kappa) = \frac{P_{\mathrm{SLNR}} D_{\kappa}}
	{\max \limits_{i} \left( \frac{1}{\sqrt{\beta_i}}\| e_i^* P_{\mathrm{SLNR}}
		D_{\kappa} \|_2 \right)
	} \label{eq:pwr2}
\end{equation}
as a reference surface parameterized by nonnegative $\kappa$ that sum into one. $P_{\mathrm{SLNR}}$ maintains a complexity of $O(m_{tx} m_{ue}^2)$ and collapses into $P_{\mathrm{ZF}}$ when $\omega \rightarrow 0$, but commonly offers superior SINR values for nontrivial $\omega$.

\item {\it Quality-of-Service convex solver}. (\cite{qos1, bjornsson, qos2, convex}) For the purpose of building an ideal reference, we can select a ray direction vector $\gamma \in \mathbb{R}^{m_{ue}}_+$, and accurately intersect the ray with a Pareto boundary using a QoS solver of (\ref{eq:convexbasic}) and a line-search along one positive scaling constant - which will give us one {\it correct} Pareto-boundary point. This reference is very slow\footnote{One iteration of an interior point method for (\ref{eq:convexbasic}) has a naive complexity of $O(m_{tx}^3 m_{ue}^3)$, reducible down to fourth order complexity using the innate matrix structures.} and unusable in real-time, but we will use it to verify that our developed parametric precoding scheme accurately arrives at the true Pareto boundary.
\end{itemize}

\begin{figure*}[h]
	\centering
	\subfloat[$H \in \mathbb{C}^{60 \times 24}$]{\includegraphics[width=0.45\linewidth]{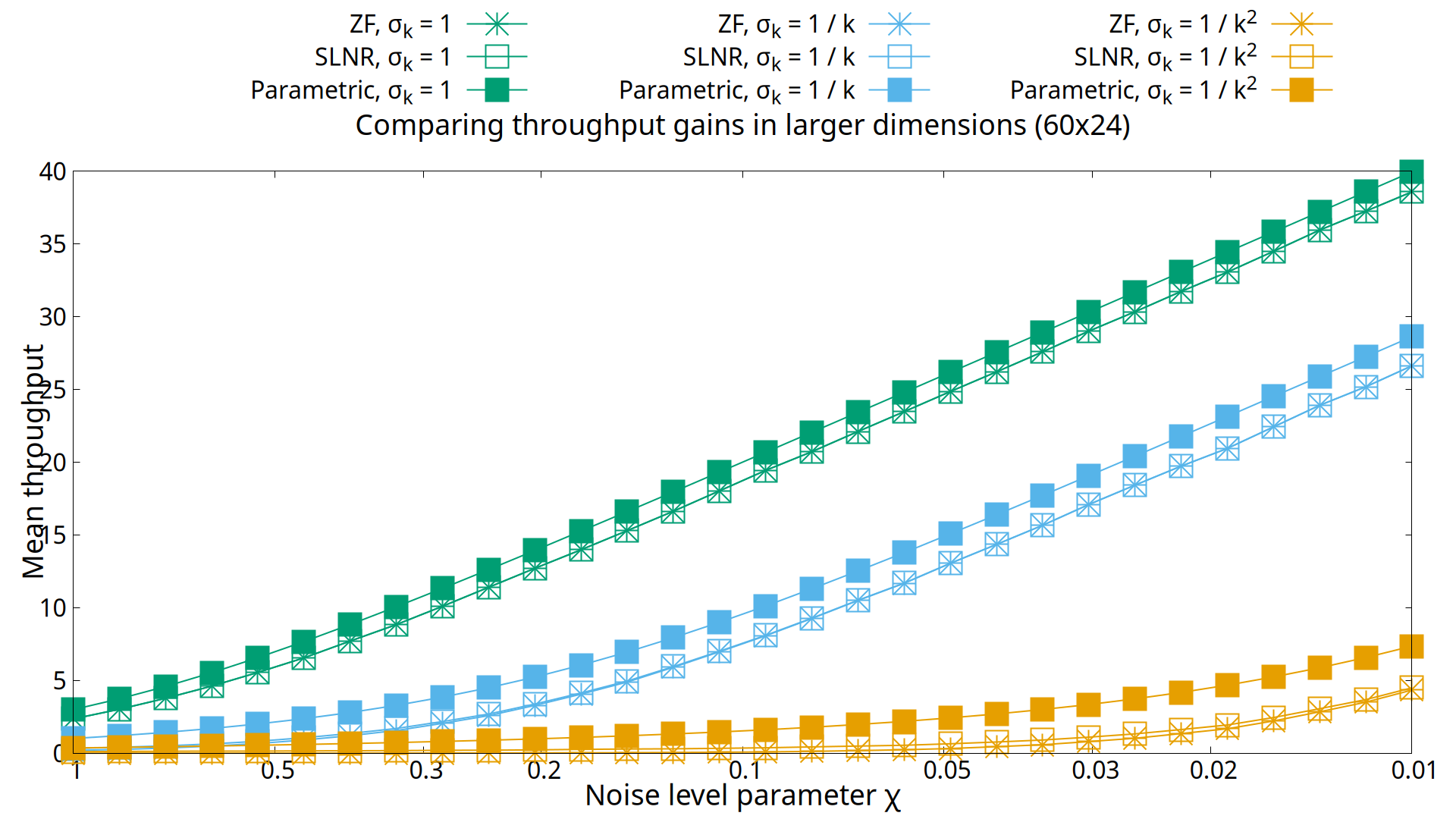}
		\label{subfig:thro_60_24}}
	\hfil
	\subfloat[$H \in \mathbb{C}^{192 \times 32}$]{\includegraphics[width=0.45\linewidth]{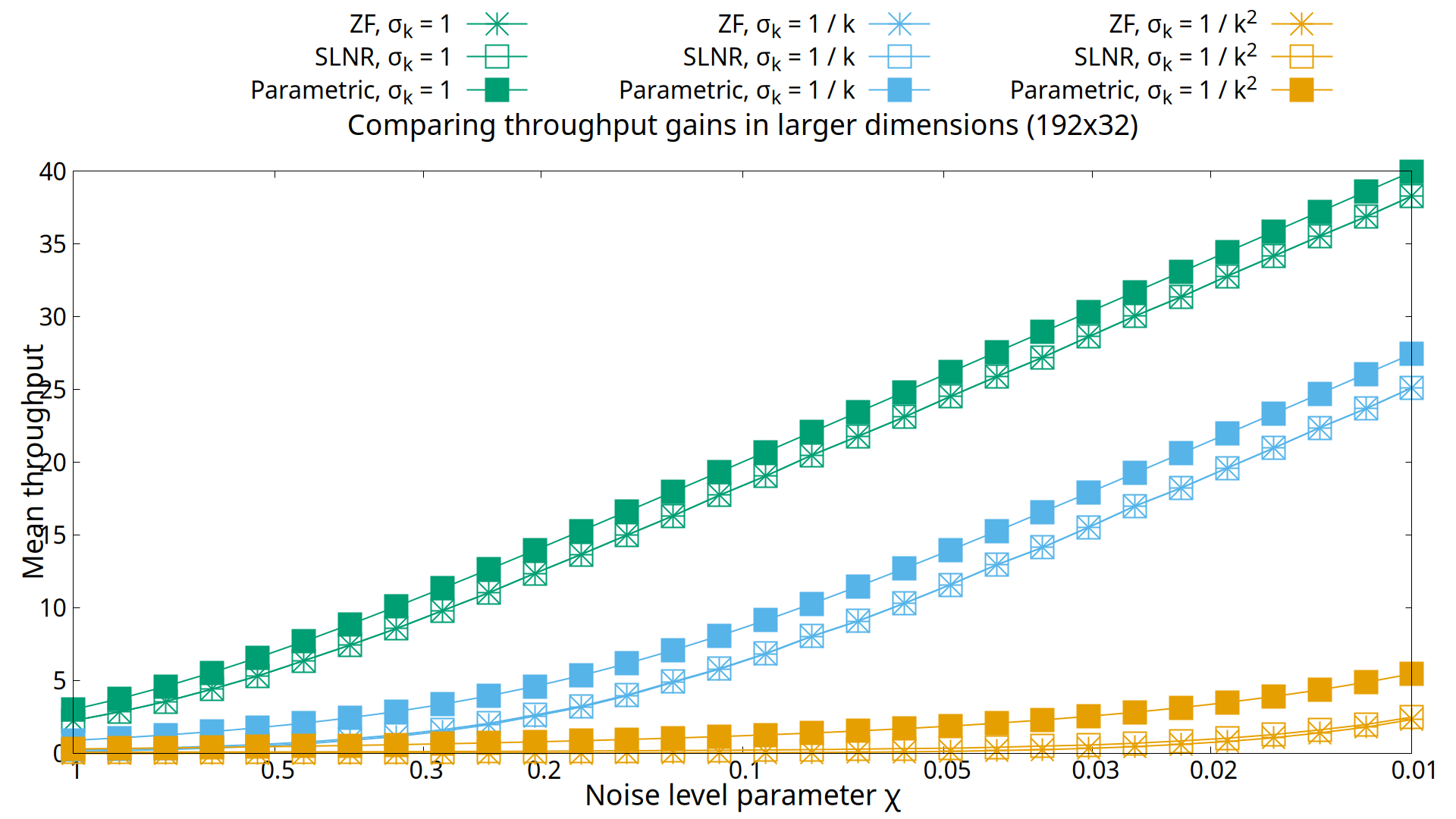}
		\label{subfig:thro_192_32}}
	
	\caption{\centering Absolute mean throughput values $\frac{1}{m_{ue}}\sum \limits_{j = 1}^{m_{ue}} 10 \log_{10} (1 + \mathrm{SINR}_j)$ for varied noise level parameter $\chi$ (\ref{eq:chi}) on channels with different singular value decay laws. ZF/SLNR references utilize water-filling. Parametric solution utilizes $\lambda_j = \frac{1}{m_{ue}}$.}
	\label{fig:thro}
\end{figure*}

\subsection{3D experiments on a toy channel}

We will begin our numerical experiments with a toy fixed channel $H \in \mathbb{R}^{8 \times 3}$. This example is considered because it is easily reproducible and allows making explicit 3D Pareto-boundary graphs - if the reader is only interested in simulations with complex channels and higher user numbers, we refer to Subsection \ref{subsec:asymp} and further.

The Figure \ref{fig:toy} provides the considered parameters and the surfaces obtained by the proposed algorithm and the two references of ZF and SLNR. In order to show that a large number of iterations is not really necessary for the proposed algorithm in applications - we are providing three surfaces:
\begin{itemize}
	\item A 'starter' surface where the parametric precoder was only computed for starter Lagrange weights $\mu = \frac{1}{m_{tx}}$ with no optimization of $\mu$ - hence each point is generated by running Algorithm \ref{alg:twopar} just once;
	\item A 'one-iteration' surface where the $\mu$ weights are updated once according to (\ref{eq:muupdate}), and hence each point is generated by running Algorithm \ref{alg:twopar} twice;
	\item A 'convergent' surface, where the $\mu$ weights are updated till the convergence condition (\ref{eq:muexit}) holds with $\delta = 0.01$.
\end{itemize}
It can be seen that running even one $\mu$ update iteration gives a quite competitive result, hence the number of iterations used can be leveraged to balance between optimality and complexity easily in applications. Each of the surfaces provided is plotted pointwise, by using $1.25 \times 10^7$ points. The arrow on each graph shows {\it a discrete point} at which the highest value of the mean throughput was encountered. If the reader wants to reproduce the results, we are providing detailed precoding matrix values for these points in Supplementary materials.

\begin{figure*}[!ht]
	\centering
	\subfloat[$H \in \mathbb{C}^{60 \times 24}$]{\includegraphics[width=0.45\linewidth]{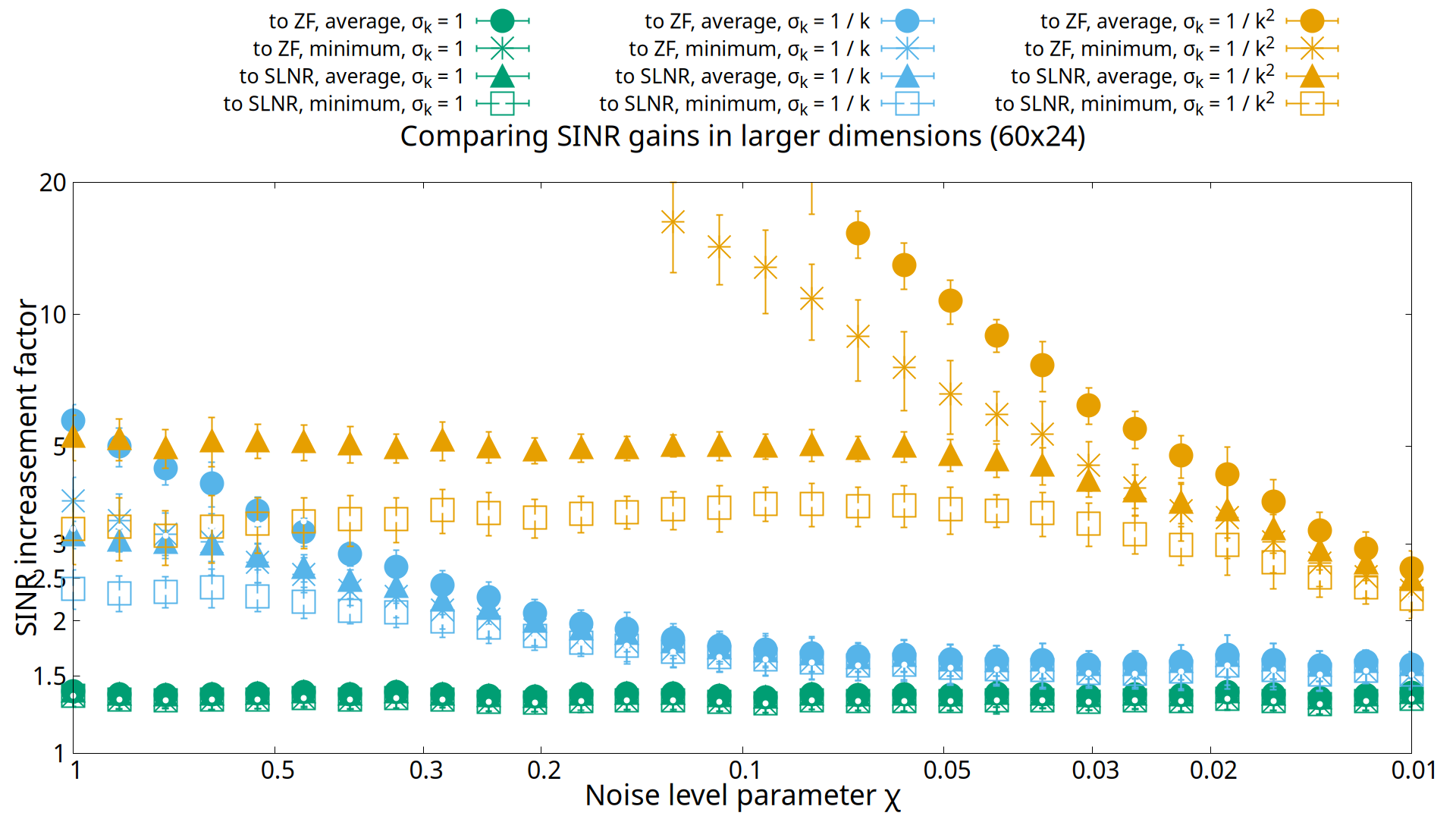}
		\label{subfig:varn_60_24}}
	\hfil
	\subfloat[$H \in \mathbb{C}^{192 \times 32}$]{\includegraphics[width=0.45\linewidth]{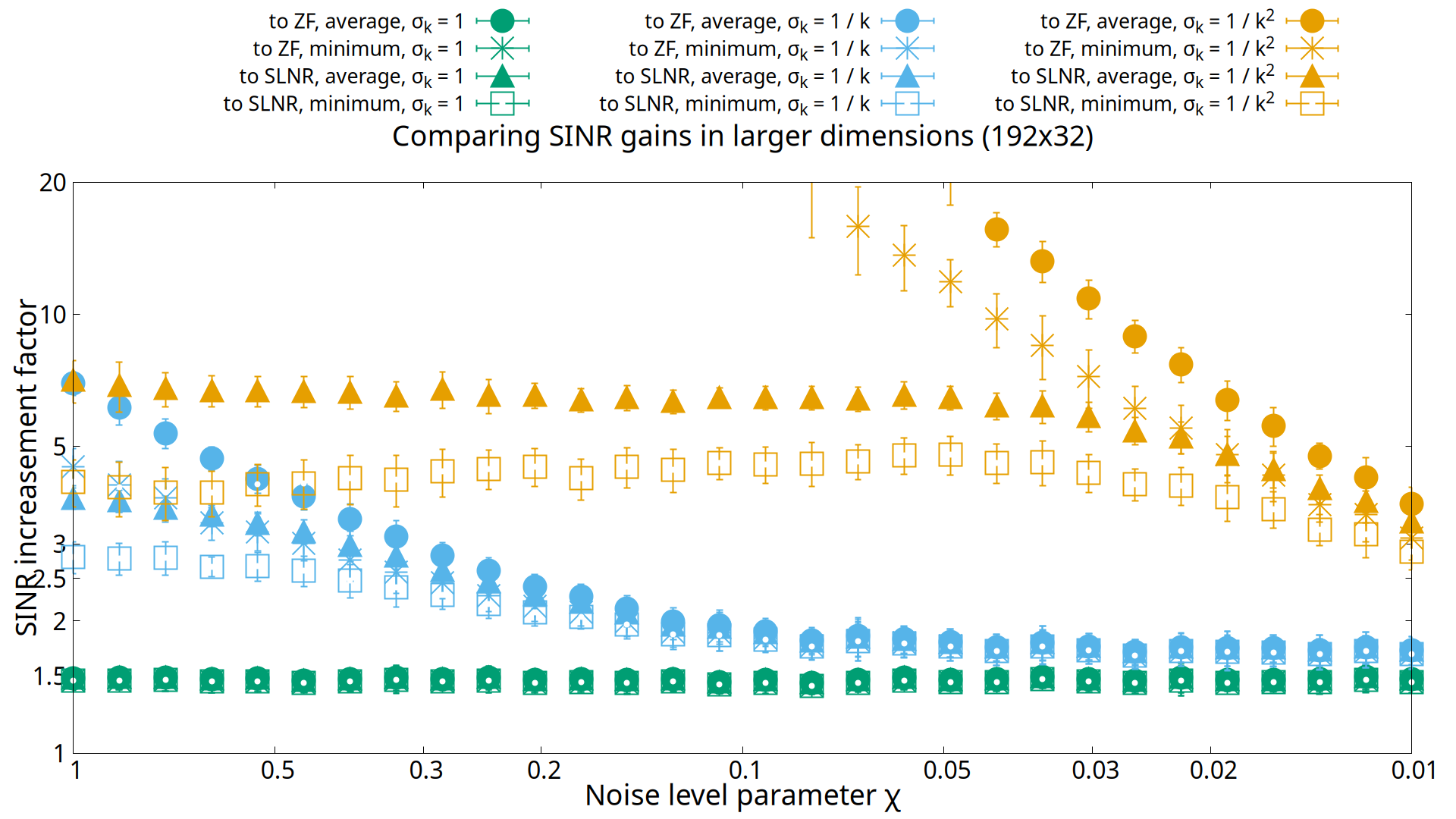}
		\label{subfig:varn_192_32}}
	
	\caption{\centering Per-SINR relative gains ${\cal G}_{avg}, {\cal G}_{min}$ (\ref{eq:usergain}) for varied noise level parameter $\chi$ (\ref{eq:chi}) on channels with different singular value decay laws. ZF/SLNR references utilize uniform power allocation. Parametric solution utilizes $\lambda_j = \frac{1}{m_{ue}}$.}
	\label{fig:varn}
\end{figure*}

\subsection{Pareto-optimality verification}

Let us now numerically check whether the developed algorithm does indeed parameterize the full Pareto boundary. Let us invoke the slow interior-point QoS solver coupled with a line-search as a reference and compare the resulting surfaces: such a comparison is provided in Figure \ref{fig:ip}. The surfaces on the figure are comprised using less points ($7.8 \times 10^3$), for two reasons:
\begin{itemize}
	\item Because the line search using convex solver runs is slow;
	\item In order to explain parameterization differences visually: the developed precoding algorithm depends on the Lagrange multipliers $\lambda$ in a nonlinear way. While the gereral monotonicity between $\lambda_k$ and $\mathrm{SINR}_k$ commonly holds, the exact location on the boundary cannot be determined in advance based on underlying $\lambda$.
\end{itemize}
To verify that the developed algorithm converges to Pareto-optimal points in higher dimensions, we devised the following experiment. Assume we fix the problem size and power/noise parameters $m_{tx}, m_{ue}, \omega, \beta$, but vary the channels $H \in \mathbb{C}^{m_{tx} \times m_{ue}}$ by generating elementwise i.i.d. zero-mean random Gaussian matrices. Assume we also randomly sample $\lambda$ with a uniform distribution in $[0, 1]$ and then scale those values towards a sum of 1. Assume for each sample we compute, as before, 
\begin{equation*}
P_{sample} := P(\lambda_{sample}, \mu(\lambda_{sample}))
\end{equation*}
approximately by doing iterative runs of Algorithm \ref{alg:twopar}, with $\mu$ updated according to Subsection \ref{subsec:refine}, until the refinement stopping criterion (\ref{eq:muexit}) holds with a predefined $\delta$. Then, we compute the resulting SINR values
\begin{equation*}
\mathrm{SINR}_1(P_{sample}), \dots, \mathrm{SINR}_{m_{ue}}(P_{sample})
\end{equation*}
and run the Quality-of-Service solver with values
\begin{equation*}
	(1 + \upsilon)\mathrm{SINR}_1(P_{sample}), \dots, (1 + \upsilon)\mathrm{SINR}_{m_{ue}}(P_{sample})
\end{equation*}
for small $\upsilon > 0$ that is comparable to $\delta$ - in order to estimate how close $P_{sample}$ is to the achievable boundary. The results are provided in Table \ref{tab:improve}. It can be seen that the algorithm tolerance $\delta$ reflects how close the generated precoder is to the Pareto boundary quite accurately.

\subsection{Varying ${\cal N} = \omega^2$ values: correctness of Proposition \ref{th:power}}
%The condition (\ref{eq:thcond}) of the proved full power necessity Proposition \ref{th:power} can be seen as a lower bound on noise values ${\cal N}_k = \omega^2$ if the remaining problem parameters (channel $H$, power constraints $\beta$) are fixed\footnote{Note that varying noise variances $\omega_k^2$ is equivalent to varying the power constraint values $\beta_k$ by the definition of SINR as in (\ref{eq:sinr}).}. In previous experiments, $\omega$ was set to a high value of $1$, and experimentally every single run of Algorithm \ref{alg:twopar}, coupled with $\mu$-refinement of Subsection \ref{subsec:refine}, has converged.

\begin{figure*}[ht]
	\centering
	\subfloat[Parametric precoder, higher $\omega = 1$]{\includegraphics[width=0.3\linewidth]{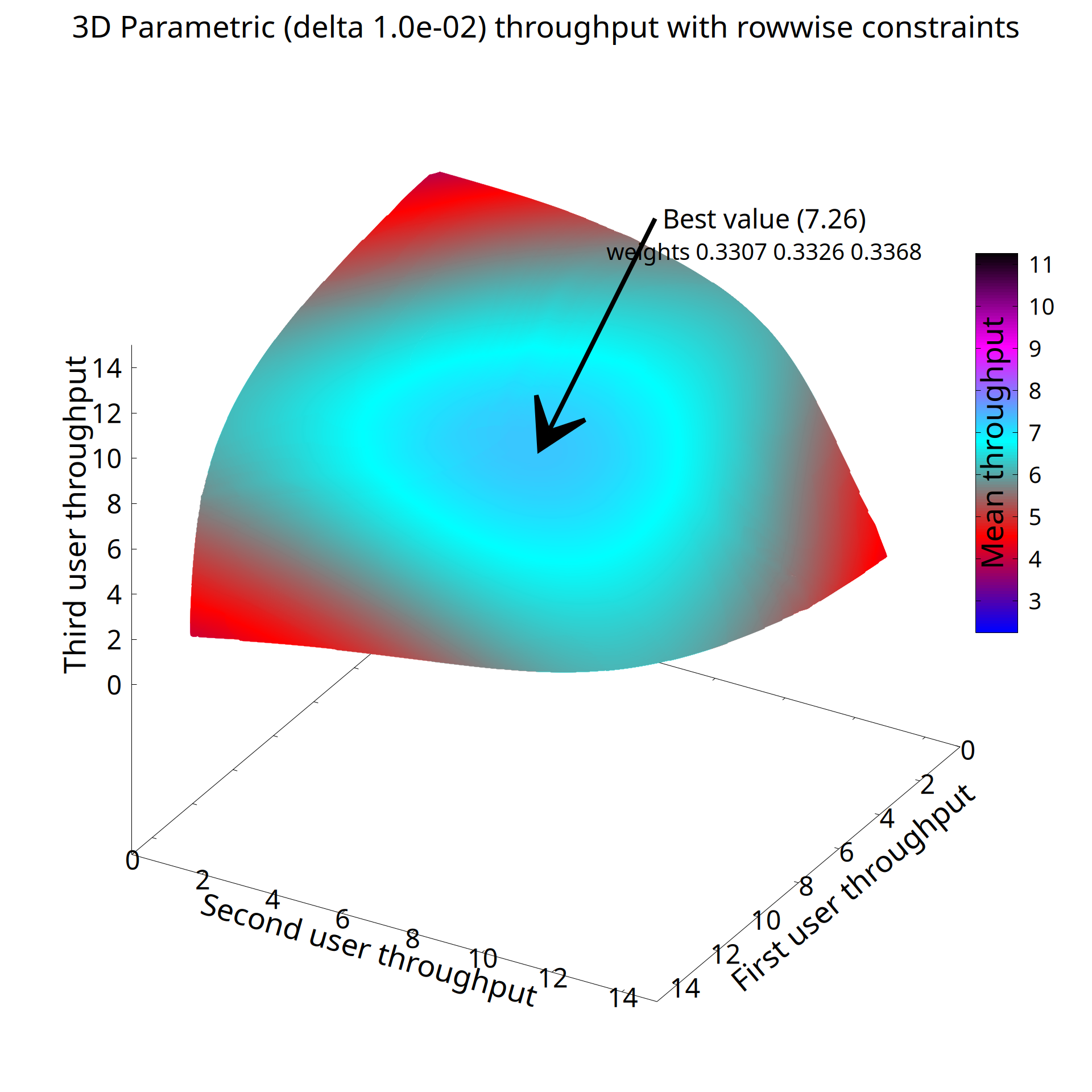}
		\label{subfig:par_high}}
	\hfil
	\subfloat[Parametric precoder, medium $\omega = 0.2$]{\includegraphics[width=0.3\linewidth]{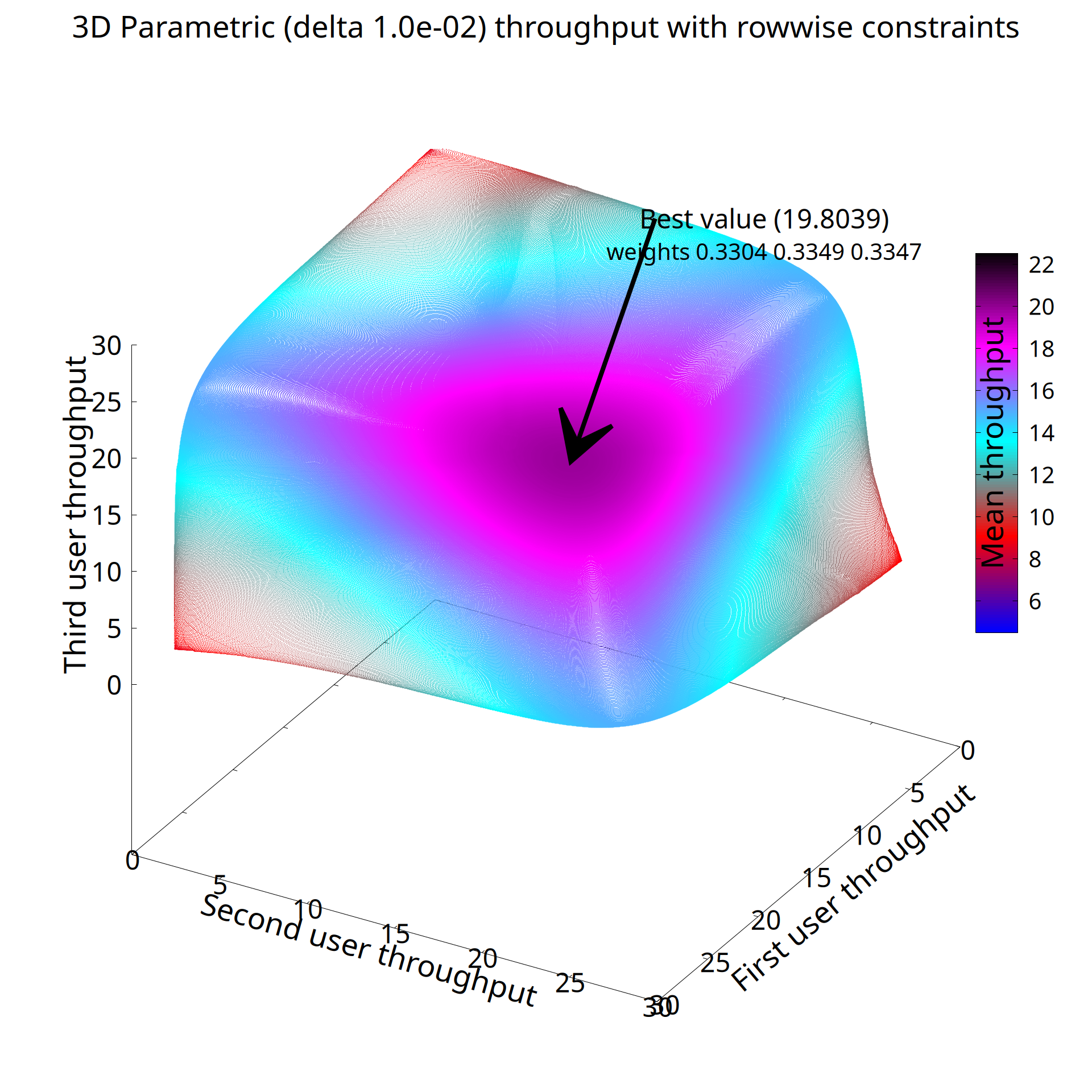}
		\label{subfig:par_medium}}
	\hfil
	\subfloat[Parametric precoder, lower $\omega = 0.04$ ]{\includegraphics[width=0.3\linewidth]{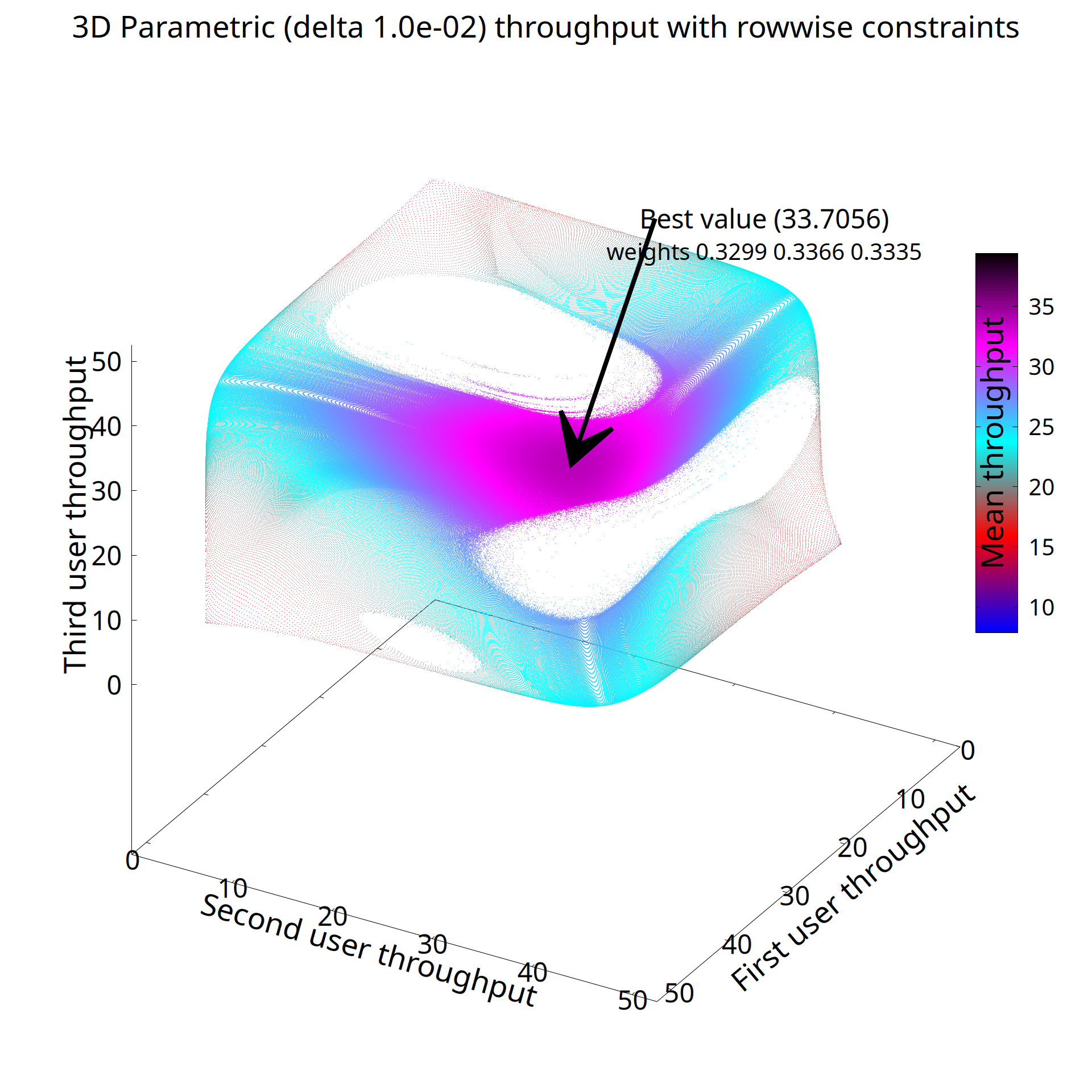}
		\label{subfig:par_low}}
	
	\centering
	\subfloat[SLNR reference, higher $\omega = 1$]{\includegraphics[width=0.3\linewidth]{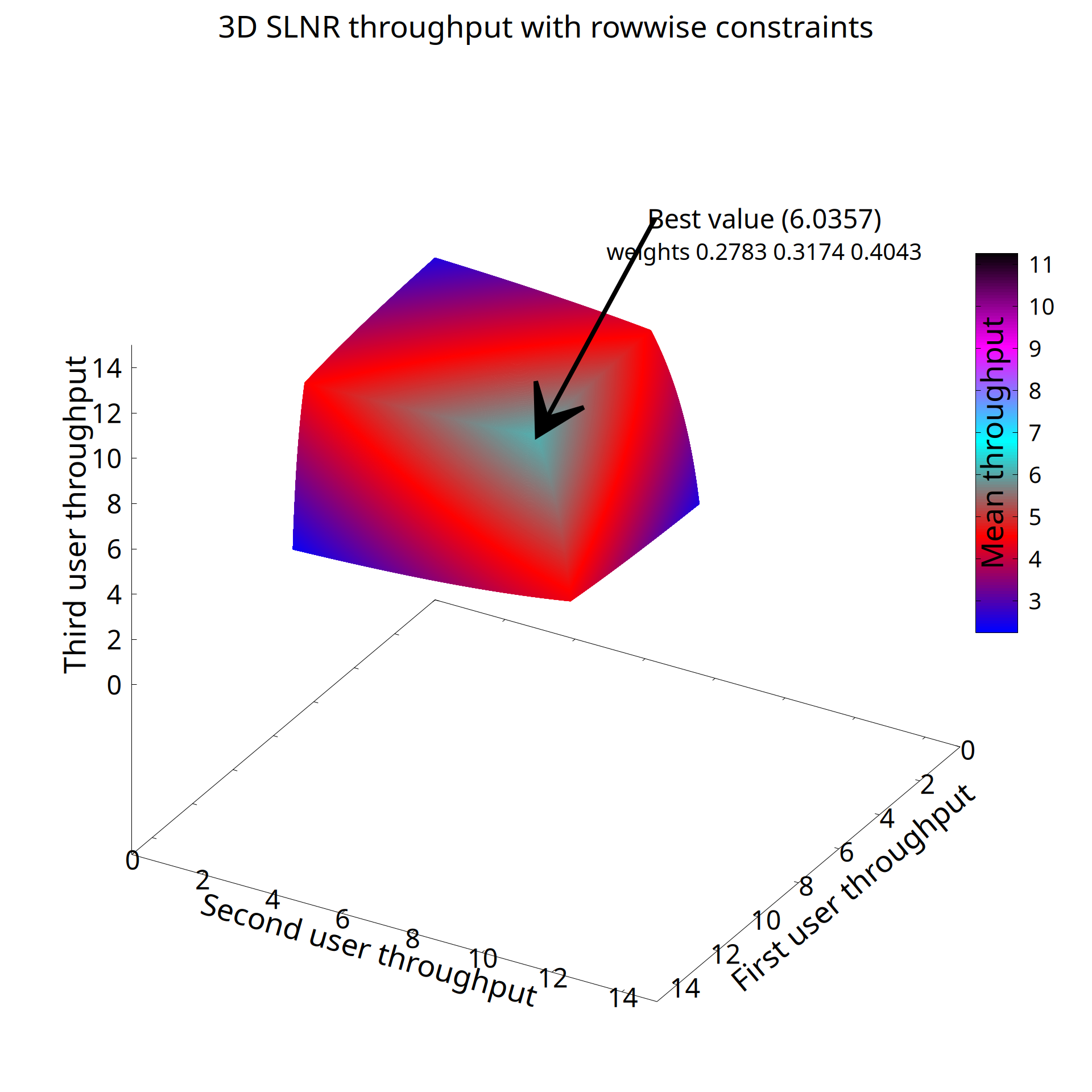}
		\label{subfig:slnr_high}}
	\hfil
	\subfloat[SLNR reference, medium $\omega = 0.2$]{\includegraphics[width=0.3\linewidth]{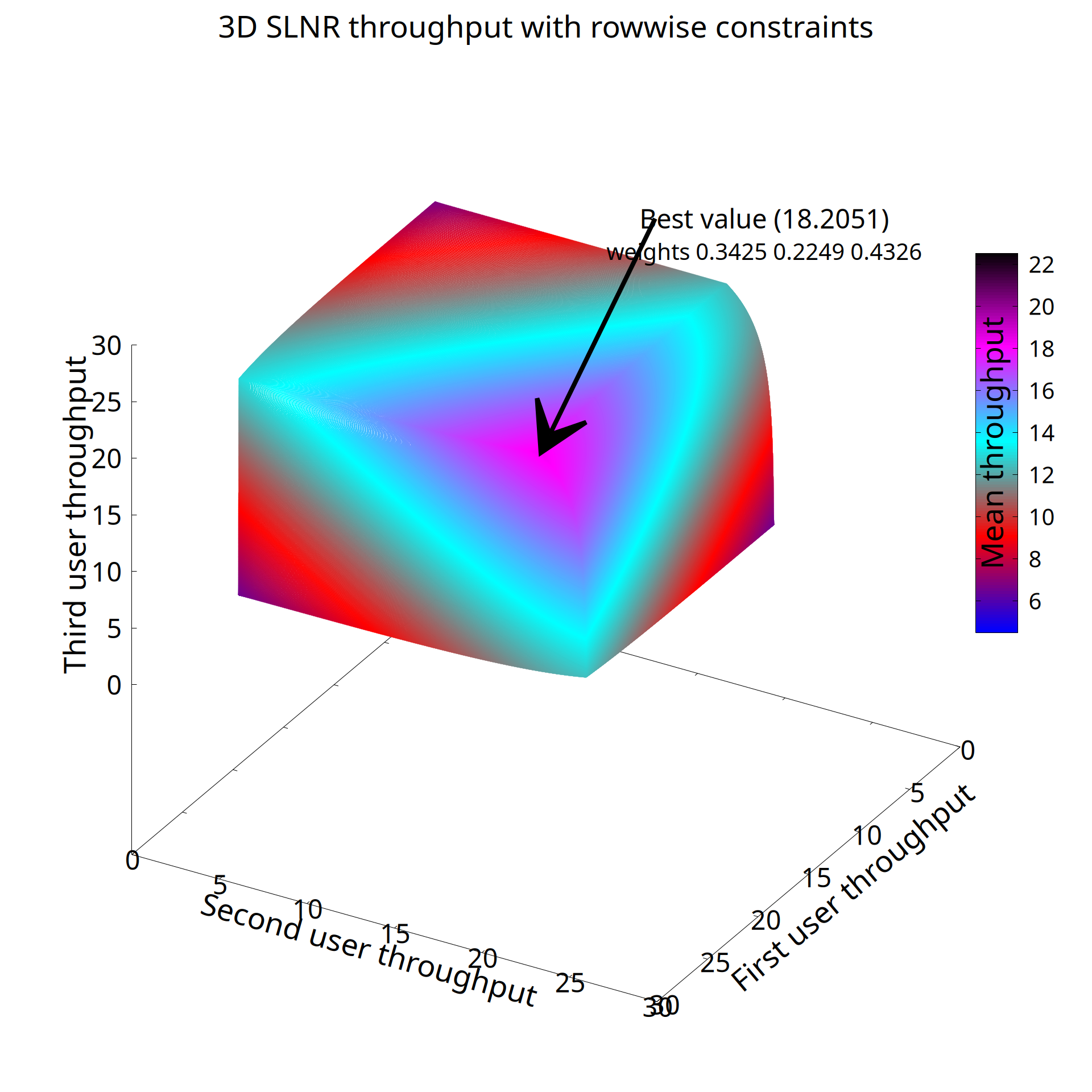}
		\label{subfig:slnr_medium}}
	\hfil
	\subfloat[SLNR reference, lower $\omega = 0.04$]{\includegraphics[width=0.3\linewidth]{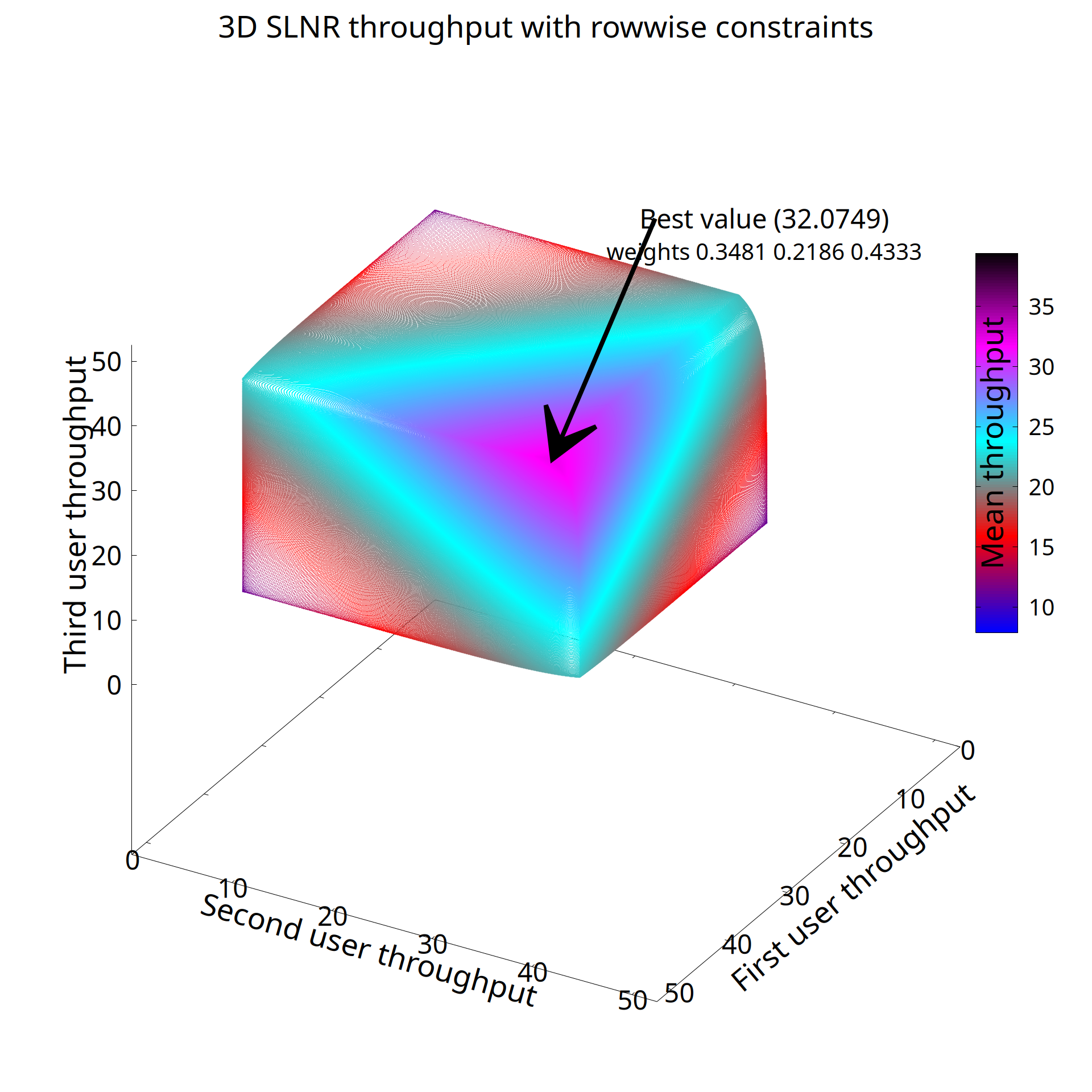}
		\label{subfig:slnr_low}}
	
	\caption{\centering Varying the SINR noise level value $\omega$ for the three-dimensional toy channel $H \in \mathbb{R}^{8 \times 3}$ given by Figure \ref{subfig:toy_ch_zf_slnr} for the developed $P(\lambda, \mu(\lambda))$ and SLNR reference. The graph slowly transforms into a cubic shape as $\omega \rightarrow 0$.}
	\label{fig:noise}
\end{figure*}

Let us study how the theory and practice break down if \mbox{$\omega \rightarrow 0$} on the same toy channel of Figure \ref{subfig:toy_ch_zf_slnr} - which would break the assumption (\ref{eq:thcond}) our theoretical analysis of Proposition \ref{th:power}. Figure \ref{fig:noise} shows how the results of the developed algorithm change for reduced $\omega$. Note that the surfaces slowly transform into a cubic shape. 

The graph that corresponds to $\omega = 0.04$ contains 'ruptures' - areas in which one of elements $\mu_i$ converges to zero. Using the 'rupture' locations we managed to extract a precoder $P_{example}$ presented in Table \ref{tab:example} in the Supplementary materials. The reader can verify that $P_{example}$ is Pareto-optimal, but doesn't utilize full power. Moreover, it was verified that Lemma \ref{lemma:unit} indeed holds and $H^* P_{example} \mathrm{diag}(v)$ has an eigenvalue \mbox{of $1$}. That confirms the main ideas of our theoretical analysis, but it doesn't mean that the developed precoding algorithm is unusable for lower noise parameters $\omega$. The computationally-efficient version of Algorithm \ref{alg:twopar} doesn't accept exact zero values of $\mu$, but if a thresholding is implemented to guarantee $\mu_j \geq \mu_{\min} > 0$, the developed algorithm still works - we did not do that on Figure~\ref{subfig:par_low} for demonstration purposes. Moreover, the most valuable points of the performance region seem unaffected by the 'ruptures' anyway.

	\begin{table*}[!t]
	\centering
	\caption{\centering Success rate of improving each SINR value by a factor $\upsilon$, 100 runs of developed algorithm with a threshold $\delta$ on random Gaussian channels}
	\label{tab:improve}
	\begin{tabular}{|c|c|c|c|c|c|}
		\hline
		\rule{0pt}{9pt} Precoder specs & $\upsilon = 10^{-3}$ & $\upsilon = 1.5 \times 10^{-4}$ & $\upsilon = 10^{-4}$ & $\upsilon = 6 \times 10^{-5}$ & $\upsilon = 10^{-5}$ \\ 
		\hline
		\rule{0pt}{9pt} $P \in \mathbb{C}^{24 \times 8}, \delta = 10^{-4}$ & 0 \% & 6 \% & 47 \% & 82 \% & 100 \% \\
		\hline
		\rule{0pt}{9pt} $P \in \mathbb{C}^{192 \times 24}, \delta = 10^{-4}$ & 0 \% &  9 \% & 28 \% & 44 \% & 95 \% \\
		\hline
	\end{tabular}
	
\end{table*}

\subsection{Varying sizes and ${\cal N} = \omega^2$ values: gain asymptotics}
\label{subsec:asymp}

Let us study the gains provided by the developed algorithm in higher dimensions $m_{tx} > m_{ue} > 3$. Let us define ground rules we will follow while comparing experimental results with different size, channel model and norm specifications.
\begin{itemize}
	\item The total Frobenius norm of the precoder will be limited to one; the respective per-antenna constraints (\ref{eq:powerrow}) will be set to equal values of $\beta_k = \frac{1}{m_{tx}}$.
	\item SINR noise coefficients $\omega$ will be selected in accordance with the channel Frobenius norm: we will utilize 
	\begin{equation}
	\label{eq:chi}
	\omega_k = \chi \frac{\|H\|_F}{m_{ue}},
	\end{equation} where $\chi$ will be a parameter.
\end{itemize}
Empirically, these rules guarantee that with a fixed $\chi$ individual throughput values $\log (1 + \mathrm{SINR}_k)$ are approximately constant regardless of $m_{tx}, m_{ue}$ - at least on randomly generated channels. The channels we will now be considering will be constructed using the SVD representation (\ref{eq:svd}) with random othogonal bases $U, V$\footnote{$U, V$ were constructed as the row/column orthogonal bases of a random Gaussian matrix in experiments.} and predefined singular value decay law in one of the three forms:
\begin{equation*}
\sigma_k = 1; \mbox{ } \sigma_k = \frac{1}{k}; \mbox{ or } \sigma_k = \frac{1}{k^2}. 
\end{equation*}

If $m_{ue} > 3$, we are unable to compare whole surfaces, and we have to restrain to comparing individual points on the surfaces shown in previous experiments. The point selection strategies will {\it roughly} optimize mean/total throughput (\ref{eq:func}):
\begin{itemize}
	\item For the developed algorithm, we will use $P_{arrow} := P(\lambda_\star, \mu(\lambda_\star))$, where ${\lambda}_{\star,k} = \frac{1}{m_{ue}}$, since it was a stronger point in most low-dimensional experiments. 
	\item For ZF and SLNR references, we will consider two strategies for 'allocating power' $\kappa$ as in (\ref{eq:pwr1}), (\ref{eq:pwr2}):
	\begin{itemize}
		\item {\it Uniform power.} The weigths $\kappa$ are selected in such a way that columns of $P_{\mathrm{ZF}} D_{\kappa}, P_{\mathrm{SLNR}} D_{\kappa}$ have equal euclidean norms.
		\item {\it Water-filling.} The weights $\kappa$ are selected in such a way\footnote{The water-filling weights that optimize ZF throughput are used for the SLNR baseline aswell, disregarding the nonzero interferences ${\cal I}_k$ of SLNR.} that $P_{\mathrm{ZF}} D_{\kappa}$ is optimal in terms of throughput (\ref{eq:func}) under Frobenius power constraint (\ref{eq:powerfro}) \cite{water}.
	\end{itemize} 
	The uniform strategy commonly offers slightly weaker throughput values, but is easier to make individual relative SINR comparisons with, because water-filling is known to set a subset of weights to exactly zero.
\end{itemize}

We proceed to analyze the values
\begin{align}
{\cal G}_{avg} := & \frac{1}{m_{ue}}\sum \limits_k \frac{\mathrm{SINR}_k(P_{arrow})}{\mathrm{SINR}_k(P_{uniform/water})} \nonumber \\
{\cal G}_{min} := & \min \limits_k \frac{\mathrm{SINR}_k(P_{arrow})}{\mathrm{SINR}_k(P_{uniform/water})}, \label{eq:usergain}
\end{align}
which are the mean and minimum gain in individual SINR vales, and the mean/total throughput values (scaled (\ref{eq:func})).  

 Let us first study the dependence of ${\cal G}_{avg}, {\cal G}_{min}$ on $\chi$ for fixed combinations $m_{tx}, m_{ue} > 3$, as provided in Figure \ref{fig:varn}. Except for poor performance of zero-forcing in high noise scenario with a poorly conditioned channel (expected), the dependence is largely constant; the constant, however, depends on the size specifications. Same conclusion can be made from the Figure \ref{fig:thro}, where the throughput values are presented: constant multiplicative gain in SINR values leads to constant additive gains in the mean throughput\footnote{If each particular throughput is increased by a constant, then the total throughput gain grows linearly with $m_{ue}$, however.}. 
  
 Let us then study the gain dependencies on the problem sizes $m_{tx}, m_{ue}$. The  Figures \ref{subfig:vars_a1}, \ref{subfig:vars_m1}, \ref{subfig:vars_a2} present example how the values (\ref{eq:usergain}) behave depending on the problem sizes $m_{tx}, m_{ue}$ for various channel singular value decay laws. We attempted to fit the experimental values with a polynomial model of the form
 \begin{equation*}
 {\cal G}_{avg}(m_{tx}, m_{ue}) \leftrightarrow c m_{tx}^{\phi_{tx}} m_{ue}^{\phi_{rx}},
 \end{equation*}
 which can be done by applying linear least squares to the value logarithms, and arrived at
 \begin{align*}
 \phi_{tx} & = 0.129, \phi_{ue} = -0.250, c = 1.781 \mbox{ for } \sigma_k := 1; \\
 	\phi_{tx} & = 0.124, \phi_{ue} = -0.209, c = 2.053 \mbox{ for } \sigma_k := \frac{1}{k}; \\
 	\phi_{tx} & = 0.119, \phi_{ue} = +0.346, c = 1.034 \mbox{ for } \sigma_k := \frac{1}{k^2}.
 \end{align*}
 Note that even the sign of $\phi_{ue}$ may change - we have to conclude that the polynomial power fittings $\phi_{tx}, \phi_{ue}$ of the gain largely depend on the channel conditioning. 

%If we proceed with even sharper singular value decays, the reference waterfilling solutions start allocating all power to very few users, the improvement becomes largely non-uniform, ${\cal G}_{min}$ drops below one, so it becomes hard to keep measuring the individual SINR gains.

\subsection{Experiments on simulated channel}
%Wireless communication algorithms are also commonly tested on software-simulated channels, designed to closer represent practical application requirements. 
We claim that no qualitative changes appear if we plug in complex-valued simulated values built using the known QuadRiga software. In order to mimic the previous experiment, we regenerated 40 channel samples using the 'BERLIN\_UMa\_NLoS' model for a large variety of size specifications\footnote{Varied $m_{ue}$ values were simulated by selecting user subsets; varied $m_{tx}$ values corresponded to separate QuadRiGa runs.}.

The results are provided on Figures \ref{subfig:vars_aq}, \ref{subfig:vars_m32}, \ref{subfig:vars_n192} along with the randomly-generated channel experiments for comparison. One can see that the SINR multiplicative gain is approximately constant (each individual SINR value commonly doubles), which leads to constant additive gain in the mean throughput values again.

\subsection{Global optimization}
While our developments {\it do not} solve the single-objective {\it global} optimization problem (\ref{eq:func}) directly, but focuses on arriving at the boundary of the feasible SINR step instead, it still helps reducing the complexity of {\it global} optimization significantly - because with the developed tool in mind it is sufficient to optimize $P(\lambda)$ over $m_{ue} - 1$ independent variables. There is a number of different strategies possible, depending on the application.

\begin{itemize}
	\item It is possible to utilize a black-box derivative-free optimizer over $\lambda$. There is a variety of options, ranging from simpler Nelder-Mead \cite{nelder} algorithm up to a sophisticated Tensor-Train optimization technique \cite{mitya}. With the Tensor-Train optimizer, we have managed to obtain a global optimum of (\ref{eq:func}) for sizes $H \in \mathbb{C}^{192 \times 24}$ in $15$ minutes of computing time on a desktop Intel Core i9-7900X processor. While this option is surely not realtime-applicable, it is significantly faster than the QoS solver coupled with the Branch-reduce-bound method, which we managed to converge with at sizes $m_{tx} < 10, m_{ue} < 10$ only.
	\item One could attempt to build a Pareto surface approximation, since in most cases this surface is either smooth or cube-like, and then find the optimal point on the approximated surface. Nonlinear interpolation models are a major trend nowadays, and some researchers \cite{deep} attempted to approximate the surface and optimize Lagrange weights previously for the case of full-power constraints (\ref{eq:powerfro}).
\end{itemize}

\section{Conclusion}

In this paper we carried out a deep study of the SINR multi-objective Pareto boundary in a special case of per-antenna power constraints. Firstly, we established a theorem that establishes that most realistic Pareto-optimal precoders utilize full power at each antenna if noise coefficients are non-negligible. Based on that, we proceeded towards building a novel iterative algorithm that converges to Pareto-optimal precoders, but has a one-iteration complexity indistinguishable from that of widely accepted Zero-Forcing/SLNR baselines. Moreover, it was shown that by altering the algorithm parameters, one can traverse the whole boundary of the feasible region, if required for various single-objective global optimization problems.

It was shown experimentally that even two iterations ('starter' step and one $\mu$ update) are sufficient to build a family of precoders that noticeably outperform ZF/SLNR baselines, while 10 iterations are sufficient to converge to the Pareto boundary accurately - so that the obtained SINR combination cannot be jointly improved by a factor of $1 + 10^{-3}$.

Our numerical study shows that the per-user relative SINR gain the Pareto-optimal precoder offer against ZF/SLNR baselines increases with $m_{tx}$ and decreases with $m_{ue}$ for well-conditioned channels, but changes mildly with the problem sizes. {\it Very roughly} that means that the per-user throughputs, in logarithmic dependence over the SINR values, are increased by a constant regardless of the initial throughput value - but this constant is larger for ill-conditioned channels.

To sum up, the developed algorithm, to our knowledge, brings the most efficient way to construct Pareto-optimal precoders {\it under per-antenna power constraints}, but still at a {\it constant-factor-higher price} as compared to the baselines of ZF/SLNR. It can be effectively used either 'as is' (by just using equal $\lambda$) or as an auxiliary tool to perform global throughput optimization over the feasible SINR region boundary. %Whether the gain offered by Pareto-optimal or globally optimal precoders with respect to baselines is worth the additional cost or not - is, of course, dependent on the application.

\begin{figure*}[!th]
		\centering
		\subfloat[Absolute mean throughput values for varied receiver number $m_{ue}$ and fixed transmitter number $m_{tx} = 192$: \mbox{$\frac{1}{m_{ue}}\sum \limits_{j = 1}^{m_{ue}} 10 \log_{10} (1 + \mathrm{SINR}_j)$}. Reference SLNR utilizes water-filling.]{\includegraphics[width=0.32\linewidth]{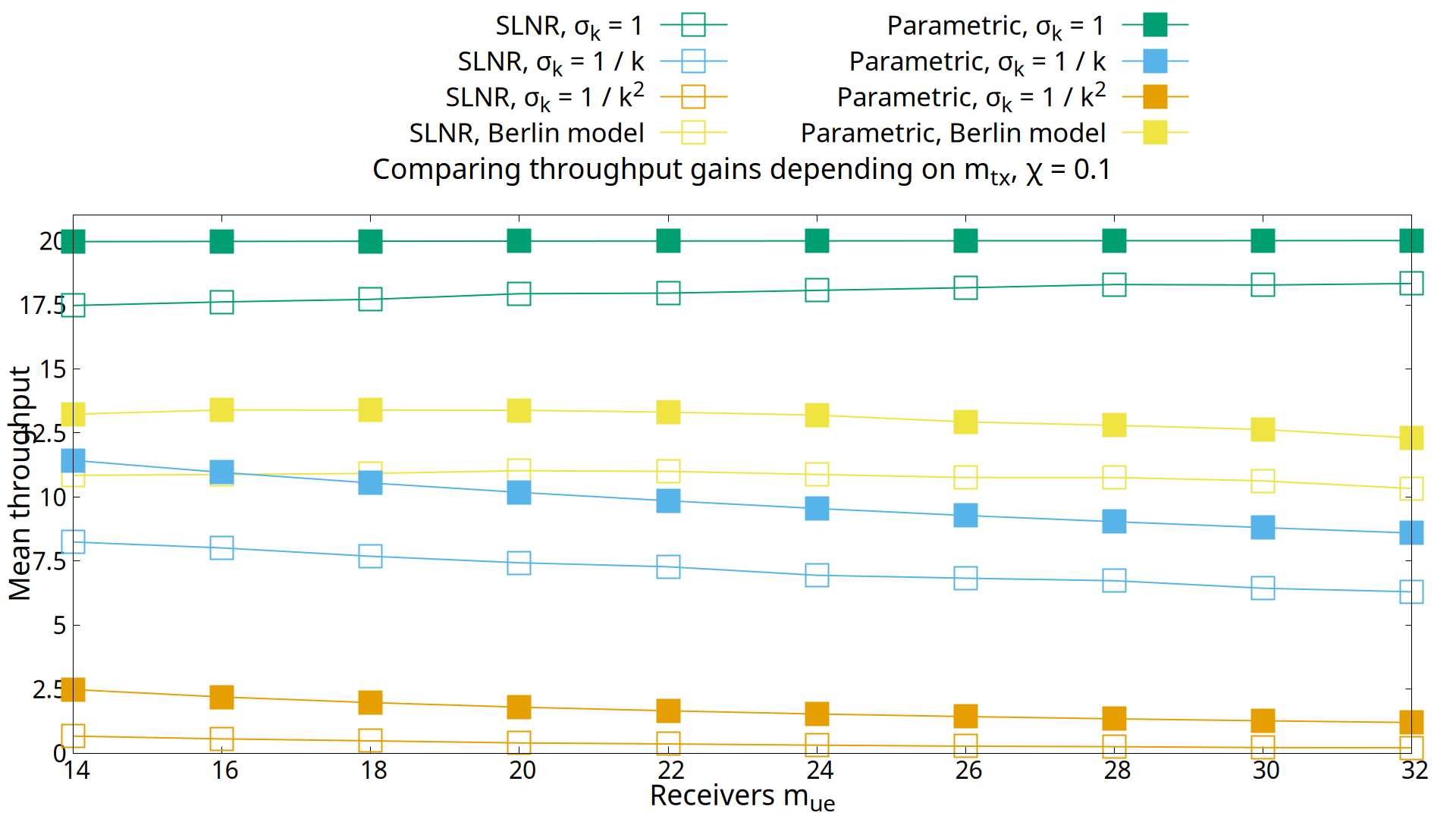}
			\label{subfig:vars_n192}}
		\hfill	
		\subfloat[Average individual SINR gain factor ${\cal G}_{avg}$ (\ref{eq:usergain}) on artificial channels with singular value decay $\sigma_k = \frac{1}{k}$. Reference SLNR utilizes uniform power allocation.]{\includegraphics[width=0.32\linewidth]{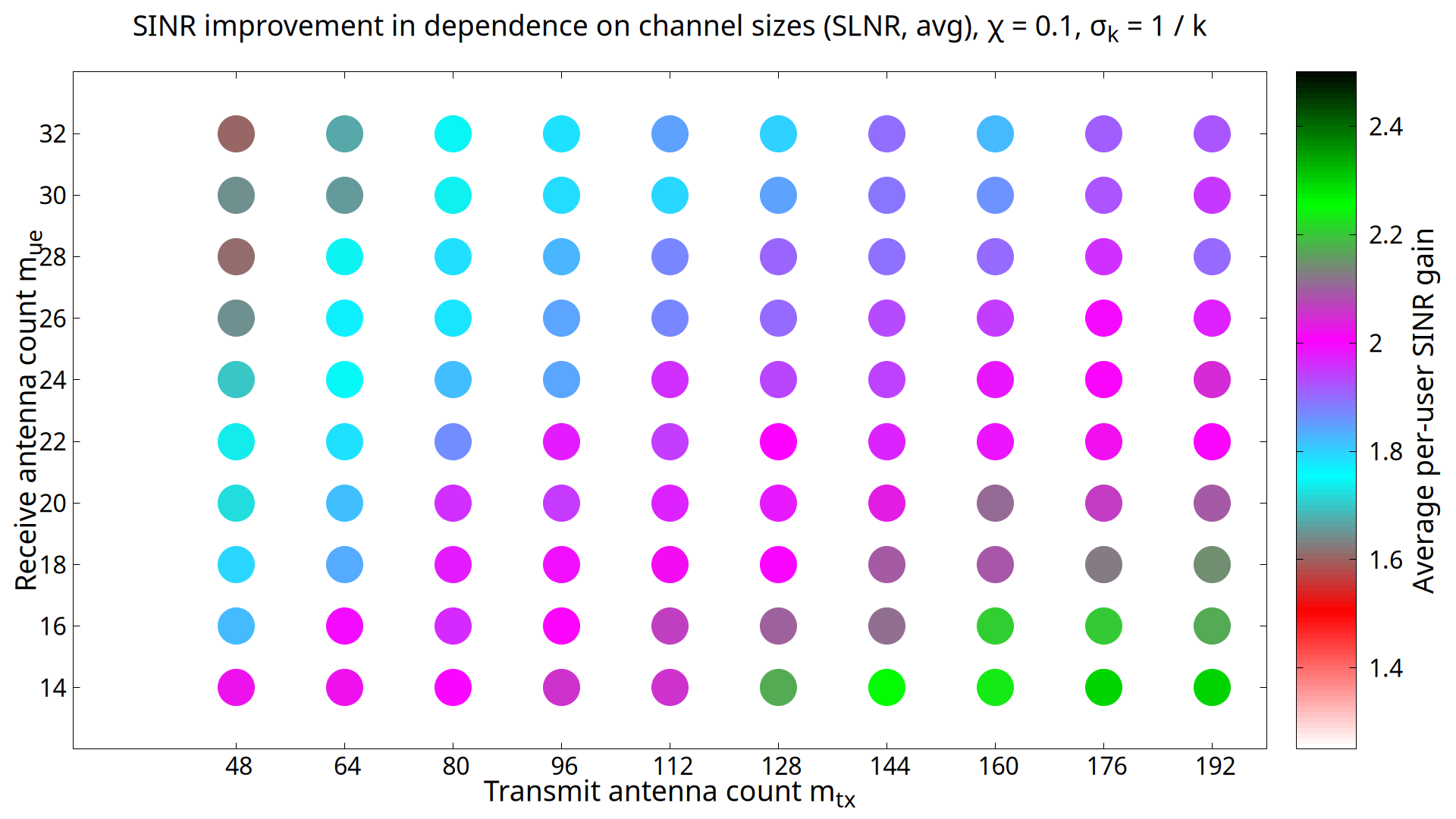}
			\label{subfig:vars_a1}}
		\hfill
		\subfloat[Minimum individual SINR gain factor ${\cal G}_{min}$ (\ref{eq:usergain}) on artificial channels with singular value decay $\sigma_k = \frac{1}{k}$. Reference SLNR utilizes uniform power allocation.]{\includegraphics[width=0.32\linewidth]{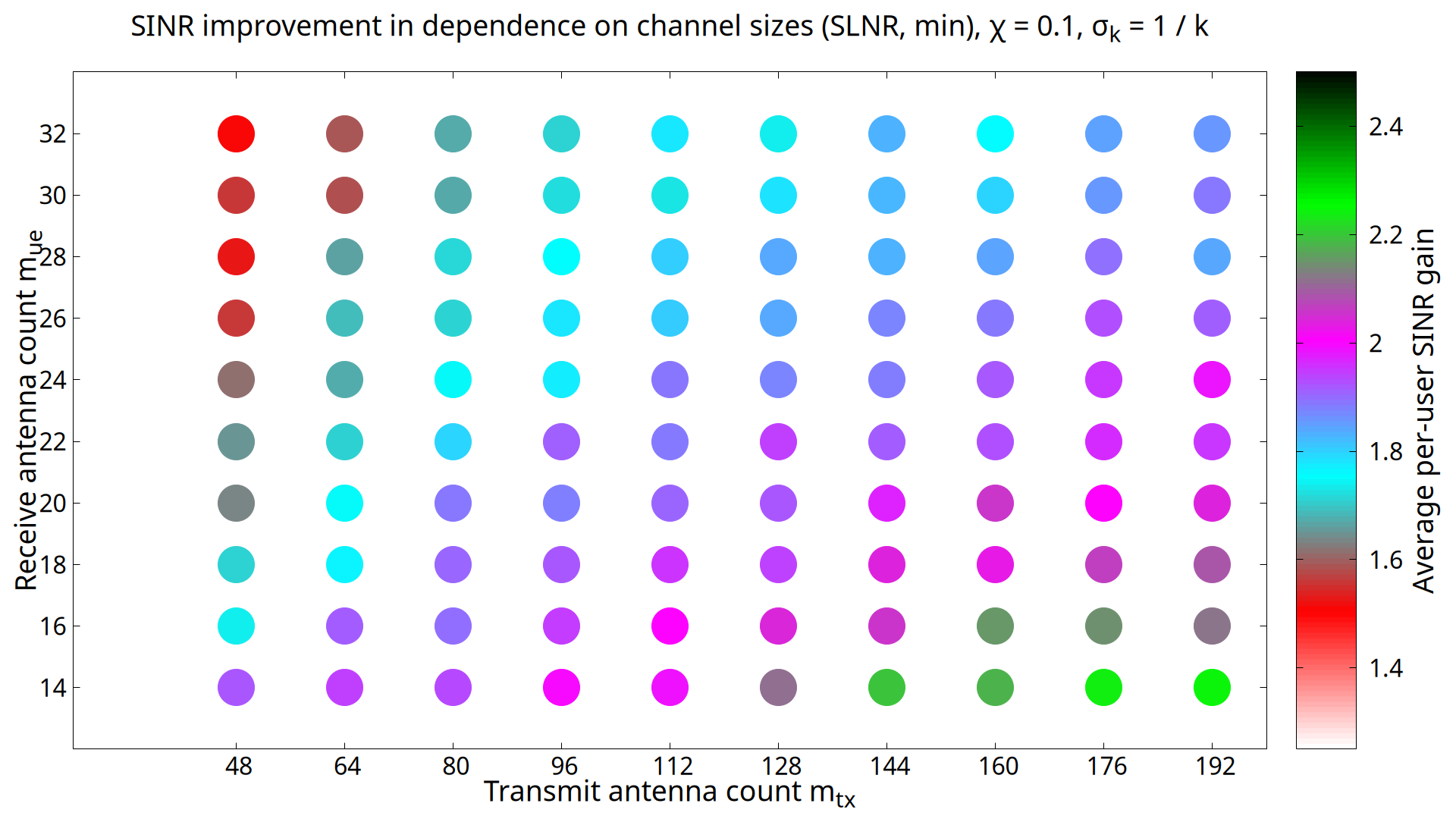}
			\label{subfig:vars_m1}}
		
			\subfloat[Absolute mean throughput values for varied transmitter number $m_{tx}$ and fixed receiver number $m_{ue} = 32$: \mbox{$\frac{1}{m_{ue}}\sum \limits_{j = 1}^{m_{ue}} 10 \log_{10} (1 + \mathrm{SINR}_j)$}. Reference SLNR utilizes water-filling.]{\includegraphics[width=0.32\linewidth]{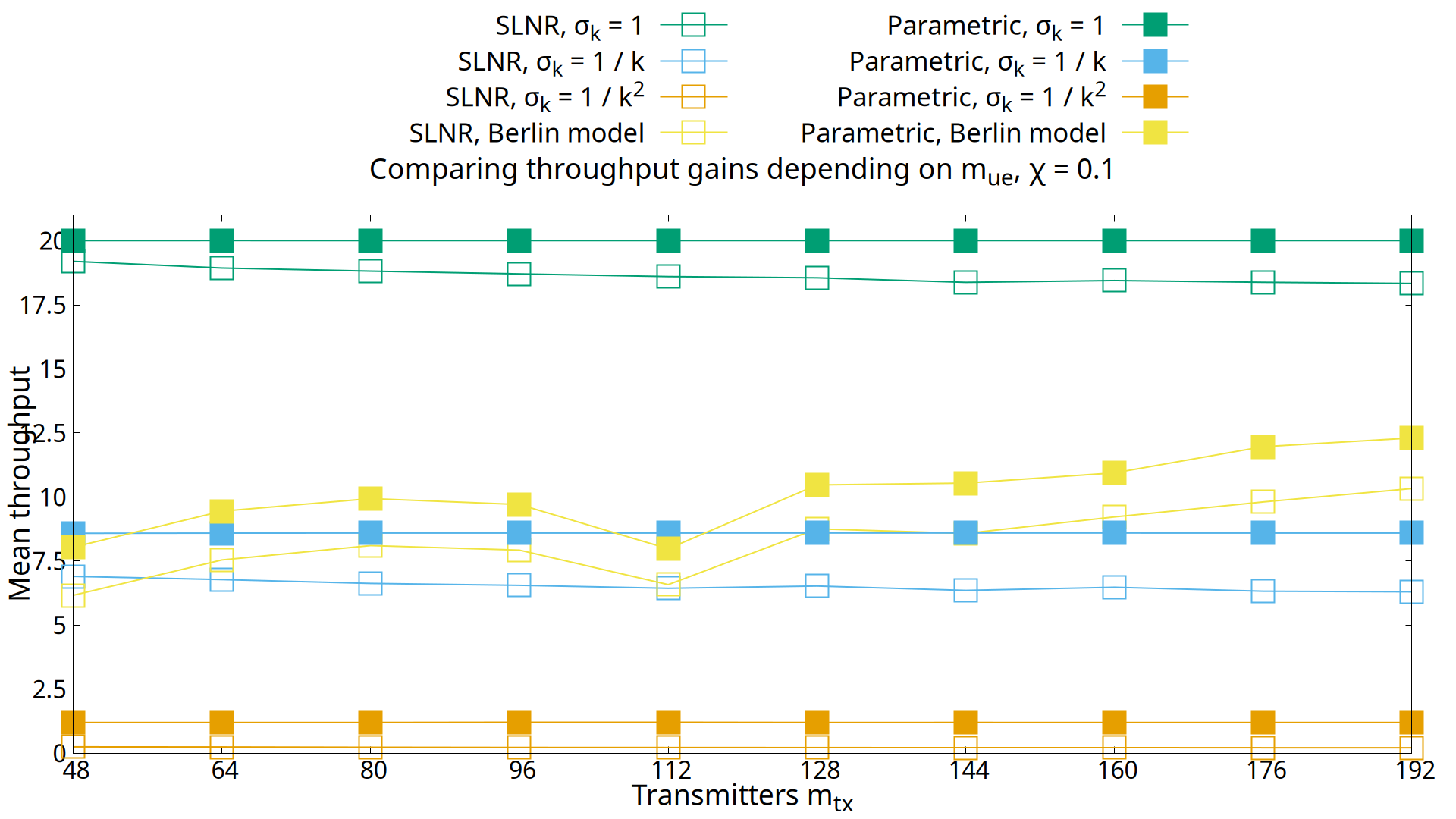}
			\label{subfig:vars_m32}}
		\hfill
		\subfloat[Average individual SINR gain factor ${\cal G}_{avg}$ (\ref{eq:usergain}) on artificial channels with singular value decay $\sigma_k = \frac{1}{k^2}$. Reference SLNR utilizes uniform power allocation.]{\includegraphics[width=0.32\linewidth]{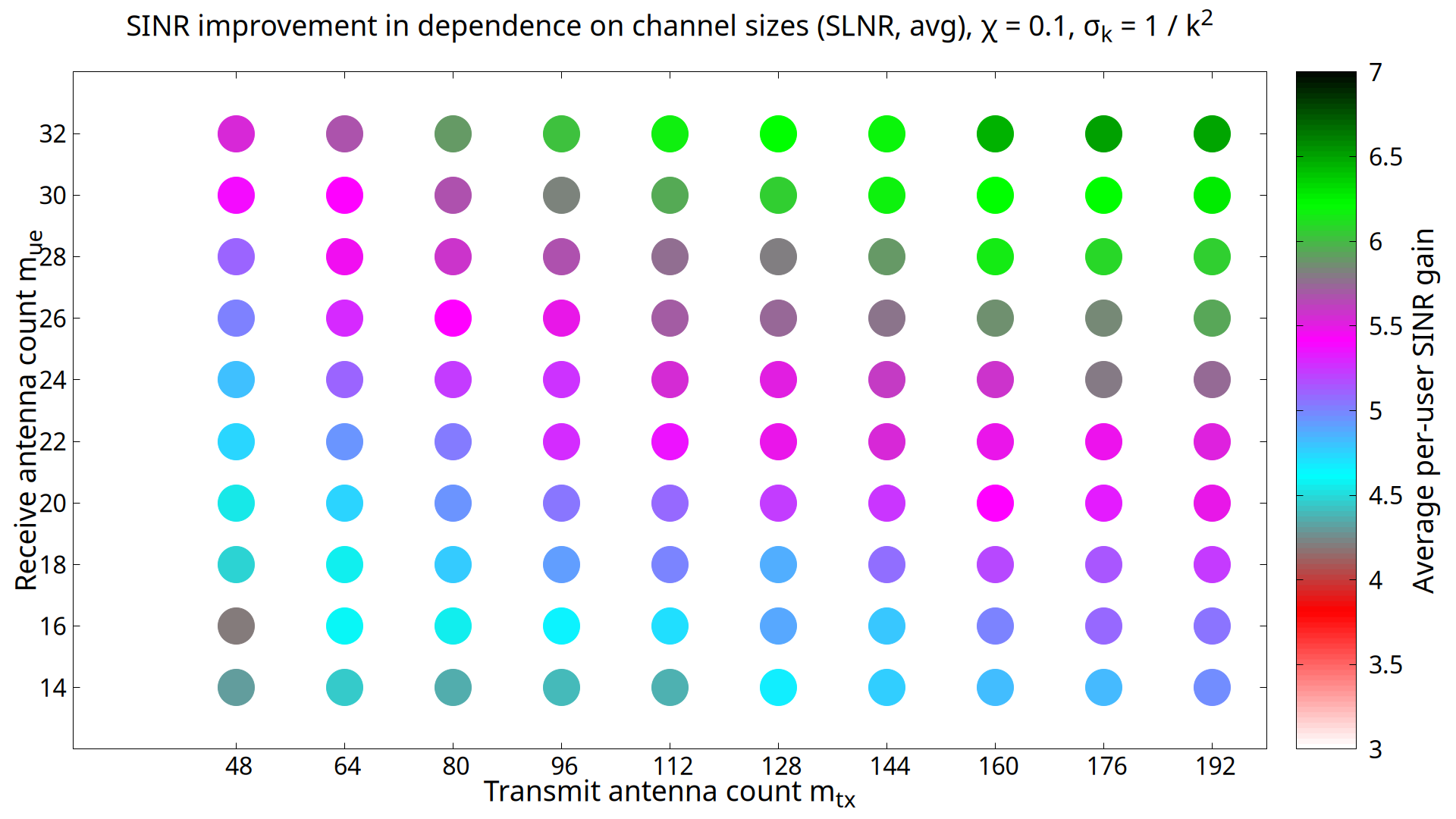}
			\label{subfig:vars_a2}}
		\hfill
		\subfloat[Average individual SINR gain factor ${\cal G}_{avg}$ (\ref{eq:usergain}) on QuadRiGa channels following the Berlin model. Reference SLNR utilizes uniform power allocation.]{\includegraphics[width=0.32\linewidth]{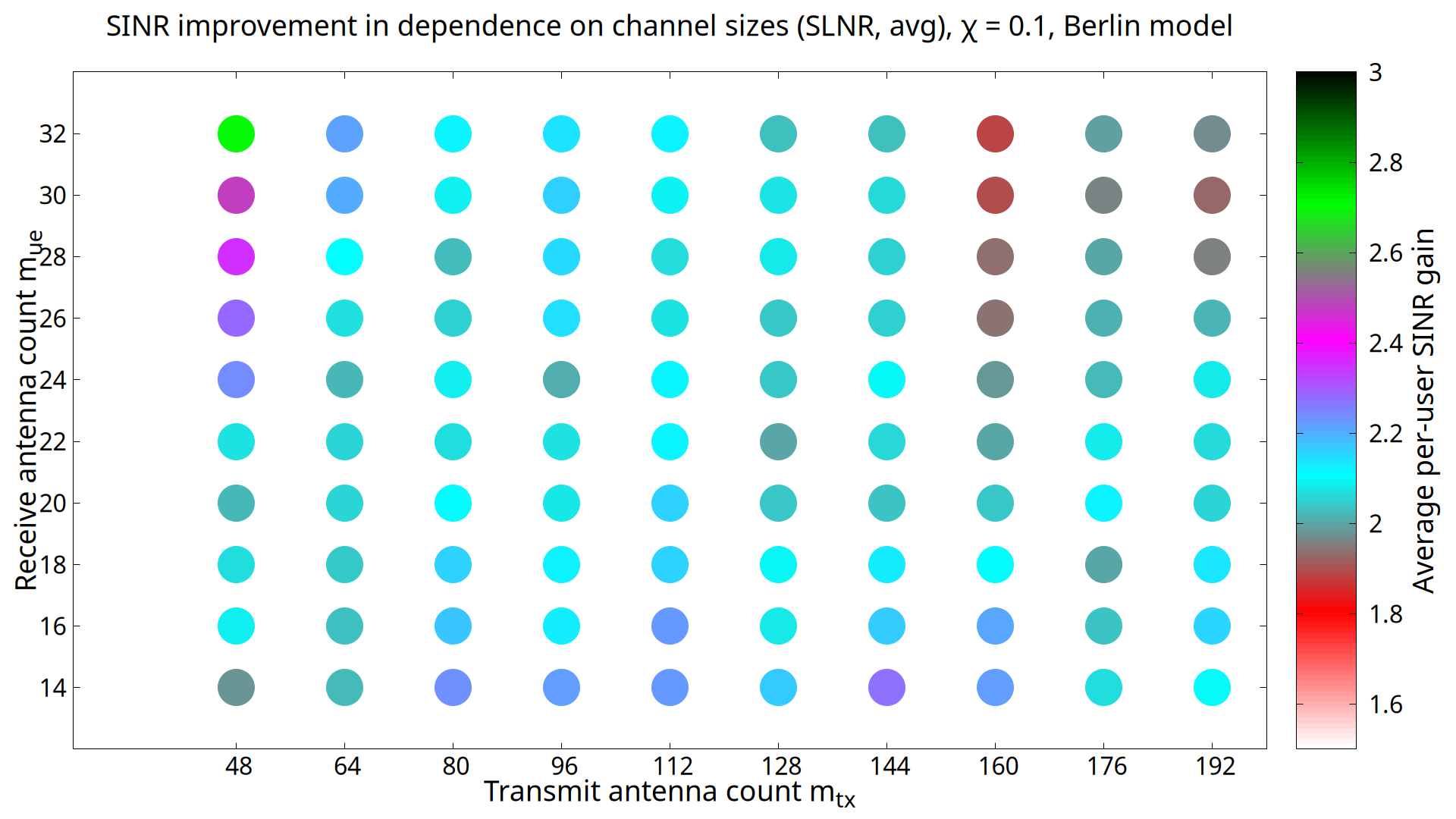}
			\label{subfig:vars_aq}}
		
		\caption{\centering Varying $m_{tx}, m_{ue}$ simultaneously with fixed SINR noise level $\chi = 0.1$ (\ref{eq:chi}) for artificial channels with various singular value decay laws and simulated channels. Each point is averaged over 40 independent experiments.
		\newline For Quadriga-generated channels, 'BERLIN\_UMa\_NLoS' model is utilized with a rectangular transmitter array with two polarizations used by each antenna. Parametric solution utilizes $\lambda_j = \frac{1}{m_{ue}}$.}
		\label{fig:vars}
	\end{figure*}

\begin{figure*}[ht]
		\centering
		\subfloat[Parametric precoder with reduced tolerance $\delta$]{\includegraphics[width=0.45\linewidth]{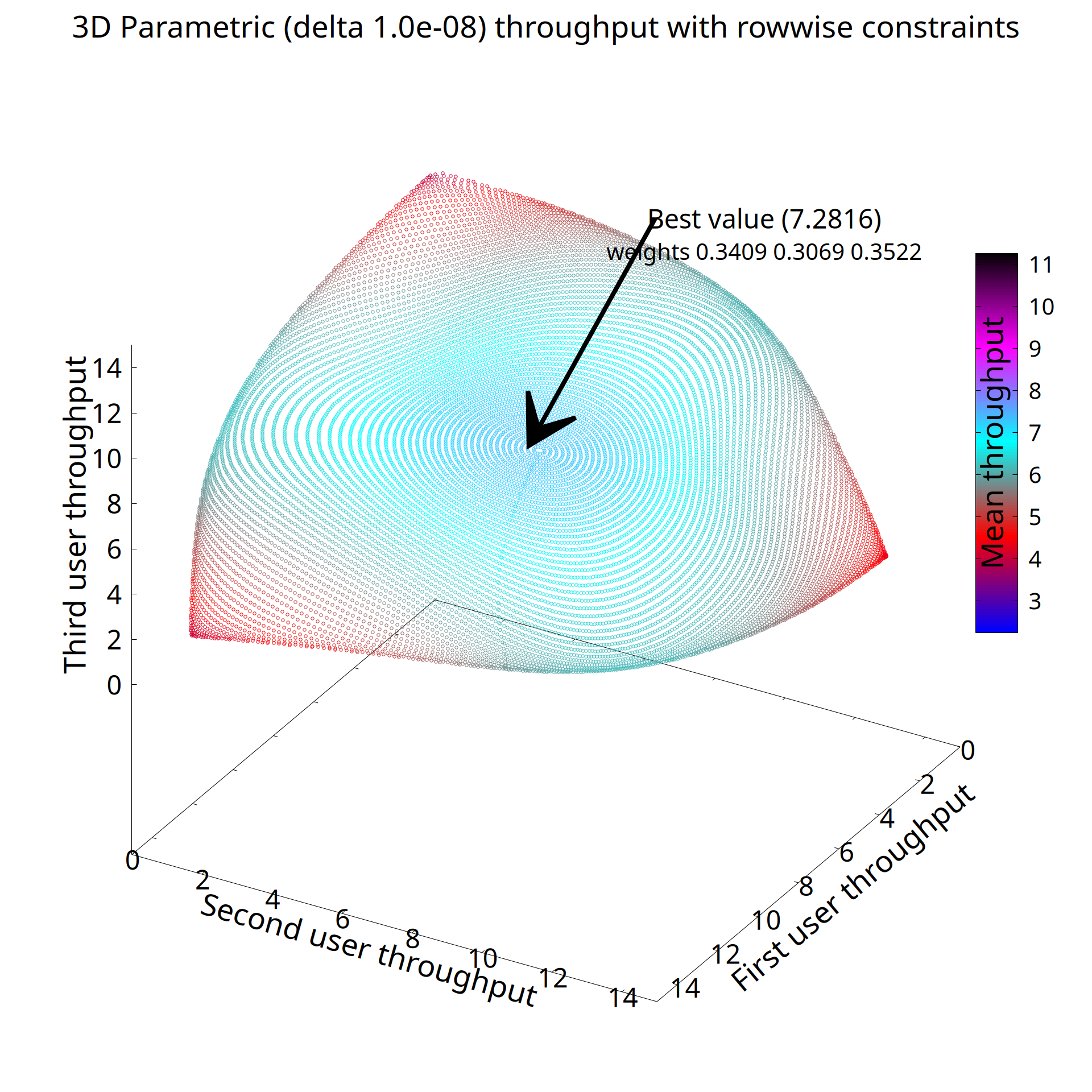}
			\label{fig_first_case}}
		\hfil
		\subfloat[Interior-point combined with line-search reference]{\includegraphics[width=0.45\linewidth]{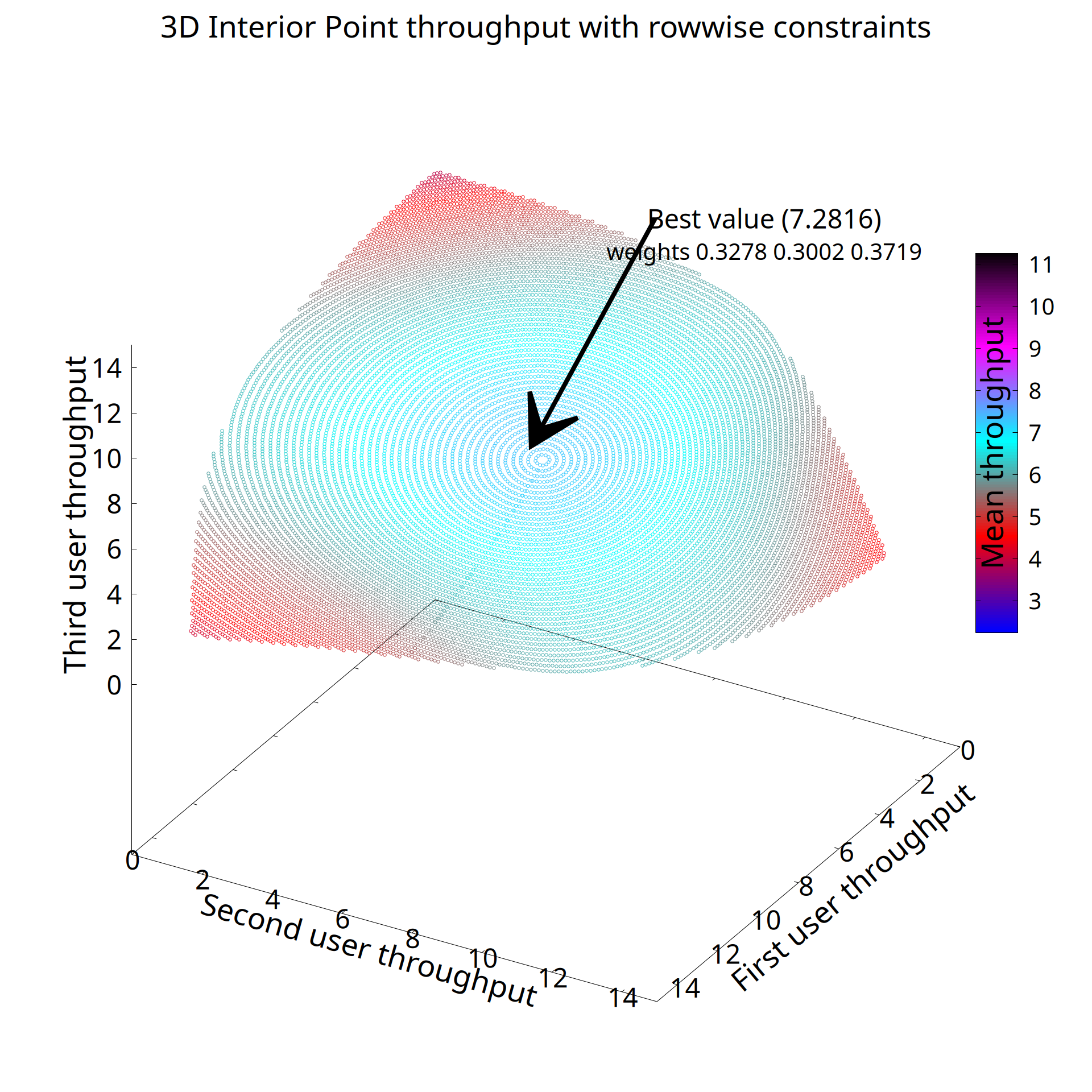}
			\label{fig_second_case}}
		
		\caption{\centering Pareto-optimality verification: comparing a surface obtained using the developed precoder $P(\lambda, \mu(\lambda))$, parameterized with $\lambda$, against an accurate Pareto boundary, parameterized with emitted line-search directions.}
		\label{fig:ip}
	\end{figure*}

{\appendix
	Due to their large double-column size, we arranged Figures \ref{fig:vars}, \ref{fig:ip} as the Appendix at the end of the paper.
}

\bibliographystyle{IEEEtran}
\bibliography{IEEEabrv, refs}

\begin{thebibliography}{10}
\providecommand{\url}[1]{#1}
\csname url@samestyle\endcsname
\providecommand{\newblock}{\relax}
\providecommand{\bibinfo}[2]{#2}
\providecommand{\BIBentrySTDinterwordspacing}{\spaceskip=0pt\relax}
\providecommand{\BIBentryALTinterwordstretchfactor}{4}
\providecommand{\BIBentryALTinterwordspacing}{\spaceskip=\fontdimen2\font plus
\BIBentryALTinterwordstretchfactor\fontdimen3\font minus
  \fontdimen4\font\relax}
\providecommand{\BIBforeignlanguage}[2]{{%
\expandafter\ifx\csname l@#1\endcsname\relax
\typeout{** WARNING: IEEEtran.bst: No hyphenation pattern has been}%
\typeout{** loaded for the language `#1'. Using the pattern for}%
\typeout{** the default language instead.}%
\else
\language=\csname l@#1\endcsname
\fi
#2}}
\providecommand{\BIBdecl}{\relax}
\BIBdecl

\bibitem{zf}
Q.~H. Spencer, A.~L. Swindlehurst, and M.~Haardt, ``Zero-forcing methods for
  downlink spatial multiplexing in multiuser mimo channels,'' \emph{IEEE
  Transactions on Signal Processing}, vol.~52, no.~2, pp. 461--471, 2004.

\bibitem{zhang1}
C.~Zhang, Y.~Jing, Y.~Huang, and L.~Yang, ``Performance analysis for massive
  mimo downlink with low complexity approximate zero-forcing precoding,''
  \emph{IEEE Transactions on Communications}, vol.~66, no.~9, pp. 3848--3864,
  2018.

\bibitem{zhang2}
J.~Liu, W.~Zhang, and Y.~Jiang, ``Fast computation of zero-forcing precoding
  for massive mimo-ofdm systems,'' \emph{IEEE Transactions on Signal
  Processing}, vol.~72, pp. 912--927, 2024.

\bibitem{mitya}
D.~Zheltkov and E.~Tyrtyshnikov, ``Global optimization based on
  tt-decomposition,'' \emph{Russian Journal of Numerical Analysis and
  Mathematical Modelling}, vol.~35, no.~4, pp. 247--261, 2020.

\bibitem{kavi}
S.~Kaviani and W.~A. Krzymie{\'n}, ``Optimal multiuser zero forcing with
  per-antenna power constraints for network mimo coordination,'' \emph{EURASIP
  Journal on Wireless Communications and Networking}, vol. 2011, no.~1, p.
  190461, 2011.

\bibitem{water}
A.~G. Armada, M.~S{\'a}nchez-Fern{\'a}ndez, and R.~Corvaja, ``Waterfilling
  schemes for zero-forcing coordinated base station transmission,'' in
  \emph{GLOBECOM 2009-2009 IEEE Global Telecommunications Conference}.\hskip
  1em plus 0.5em minus 0.4em\relax IEEE, 2009, pp. 1--5.

\bibitem{bocca}
F.~Boccardi and H.~Huang, ``Zero-forcing precoding for the mimo broadcast
  channel under per-antenna power constraints,'' in \emph{2006 IEEE 7th
  Workshop on Signal Processing Advances in Wireless Communications}.\hskip 1em
  plus 0.5em minus 0.4em\relax IEEE, 2006, pp. 1--5.

\bibitem{vu}
Q.-D. Vu, L.-N. Tran, R.~Farrell, and E.-K. Hong, ``Energy-efficient
  zero-forcing precoding design for small-cell networks,'' \emph{IEEE
  Transactions on Communications}, vol.~64, no.~2, pp. 790--804, 2015.

\bibitem{lee}
S.-R. Lee, J.-S. Kim, S.-H. Moon, H.-B. Kong, and I.~Lee, ``Zero-forcing
  beamforming in multiuser miso downlink systems under per-antenna power
  constraint and equal-rate metric,'' \emph{IEEE Transactions on Wireless
  Communications}, vol.~12, no.~1, pp. 228--236, 2012.

\bibitem{wiesel}
A.~Wiesel, Y.~C. Eldar, and S.~Shamai, ``Zero-forcing precoding and generalized
  inverses,'' \emph{IEEE Transactions on Signal Processing}, vol.~56, no.~9,
  pp. 4409--4418, 2008.

\bibitem{bjornsson}
E.~Bj{\"o}rnson, E.~Jorswieck \emph{et~al.}, ``Optimal resource allocation in
  coordinated multi-cell systems,'' \emph{Foundations and
  Trends{\textregistered} in Communications and Information Theory}, vol.~9,
  no. 2--3, pp. 113--381, 2013.

\bibitem{jors}
E.~A. Jorswieck, E.~G. Larsson, and D.~Danev, ``Complete characterization of
  the pareto boundary for the miso interference channel,'' \emph{IEEE
  Transactions on Signal Processing}, vol.~56, no.~10, pp. 5292--5296, 2008.

\bibitem{nonconv1}
Y.-F. Liu, Y.-H. Dai, and Z.-Q. Luo, ``Coordinated beamforming for miso
  interference channel: Complexity analysis and efficient algorithms,''
  \emph{IEEE Transactions on Signal Processing}, vol.~59, no.~3, pp.
  1142--1157, 2010.

\bibitem{nonconv2}
Z.-Q. Luo and S.~Zhang, ``Dynamic spectrum management: Complexity and
  duality,'' \emph{IEEE Journal of Selected Topics in Signal Processing},
  vol.~2, no.~1, pp. 57--73, 2008.

\bibitem{brb}
H.~Tuy, F.~Al-Khayyal, and P.~T. Thach, ``Monotonic optimization: Branch and
  cut methods,'' in \emph{Essays and surveys in global optimization}.\hskip 1em
  plus 0.5em minus 0.4em\relax Springer, 2005, pp. 39--78.

\bibitem{hadamard}
R.~A. Horn and Z.~Yang, ``Rank of a hadamard product,'' \emph{Linear Algebra
  and its Applications}, vol. 591, pp. 87--98, 2020.

\bibitem{woodbury}
M.~A. Woodbury, \emph{Inverting modified matrices}.\hskip 1em plus 0.5em minus
  0.4em\relax Department of Statistics, Princeton University, 1950.

\bibitem{slnr1}
M.~Sadek, A.~Tarighat, and A.~H. Sayed, ``A leakage-based precoding scheme for
  downlink multi-user mimo channels,'' \emph{IEEE Transactions on Wireless
  Communications}, vol.~6, no.~5, pp. 1711--1721, 2007.

\bibitem{slnr2}
J.~Goldberg and J.~R. Fonollosa, ``Downlink beamforming for cellular mobile
  communications,'' in \emph{1997 IEEE 47th Vehicular Technology Conference.
  Technology in Motion}, vol.~2.\hskip 1em plus 0.5em minus 0.4em\relax IEEE,
  1997, pp. 632--636.

\bibitem{qos1}
H.~Dahrouj and W.~Yu, ``Coordinated beamforming for the multicell multi-antenna
  wireless system,'' \emph{IEEE Transactions on Wireless Communications},
  vol.~9, no.~5, pp. 1748--1759, 2010.

\bibitem{qos2}
A.~Wiesel, Y.~C. Eldar, and S.~Shamai, ``Linear precoding via conic
  optimization for fixed mimo receivers,'' \emph{IEEE Transactions on Signal
  Processing}, vol.~54, no.~1, pp. 161--176, 2005.

\bibitem{convex}
S.~P. Boyd and L.~Vandenberghe, \emph{Convex optimization}.\hskip 1em plus
  0.5em minus 0.4em\relax Cambridge university press, 2004.

\bibitem{nelder}
J.~A. Nelder and R.~Mead, ``A simplex method for function minimization,''
  \emph{The Computer Journal}, vol.~7, no.~4, pp. 308--313, 1965.

\bibitem{deep}
M.~Zhang, J.~Gao, and C.~Zhong, ``A deep learning-based framework for low
  complexity multiuser mimo precoding design,'' \emph{IEEE Transactions on
  Wireless Communications}, vol.~21, no.~12, pp. 11\,193--11\,206, 2022.

\end{thebibliography}

\begin{IEEEbiography}[{\includegraphics[width=1in,height=1.25in,clip,keepaspectratio]{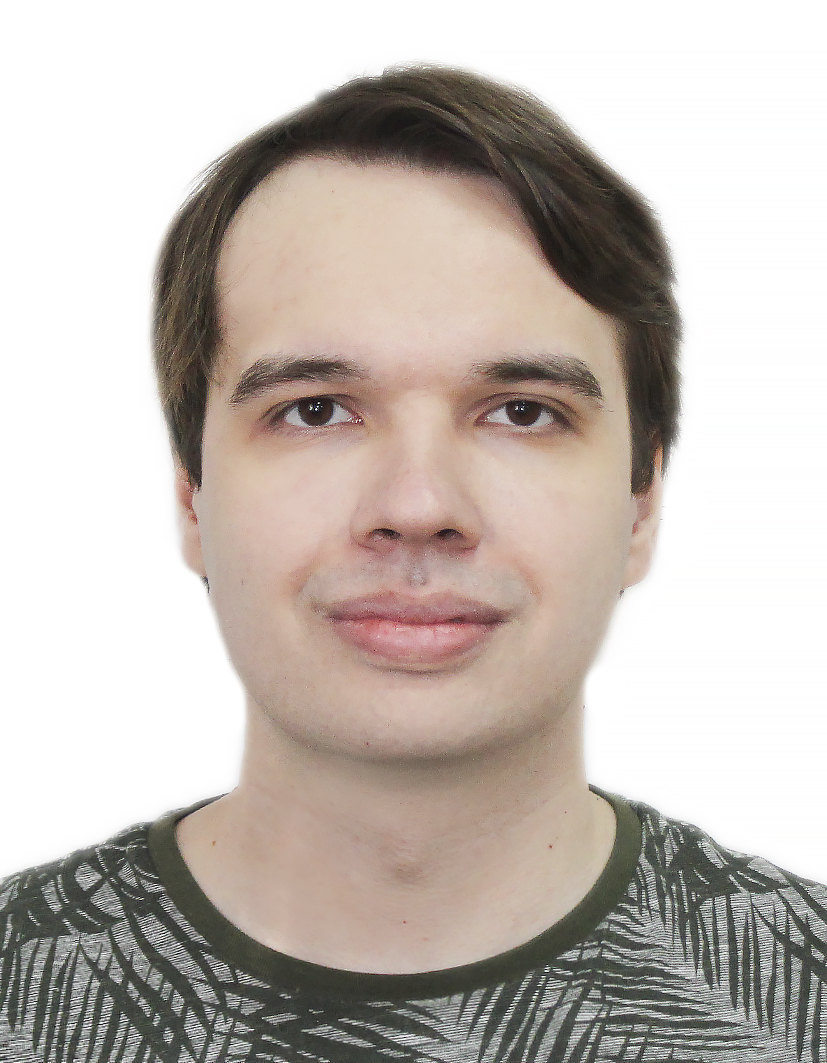}}]{Sergey Petrov} received the bachelor, master and Ph.D. degrees in Computational Mathematics in Lomonosov Moscow State University in 2016, 2018 and 2023, respectively. He is currently a post-doctoral researcher at Khalifa University, Abu-Dhabi. His research interests include computational linear algebra, optimization methods, randomized methods and perturbation theory.
\end{IEEEbiography}

\begin{IEEEbiography}[{\includegraphics[width=1in,height=1.25in,clip,keepaspectratio]{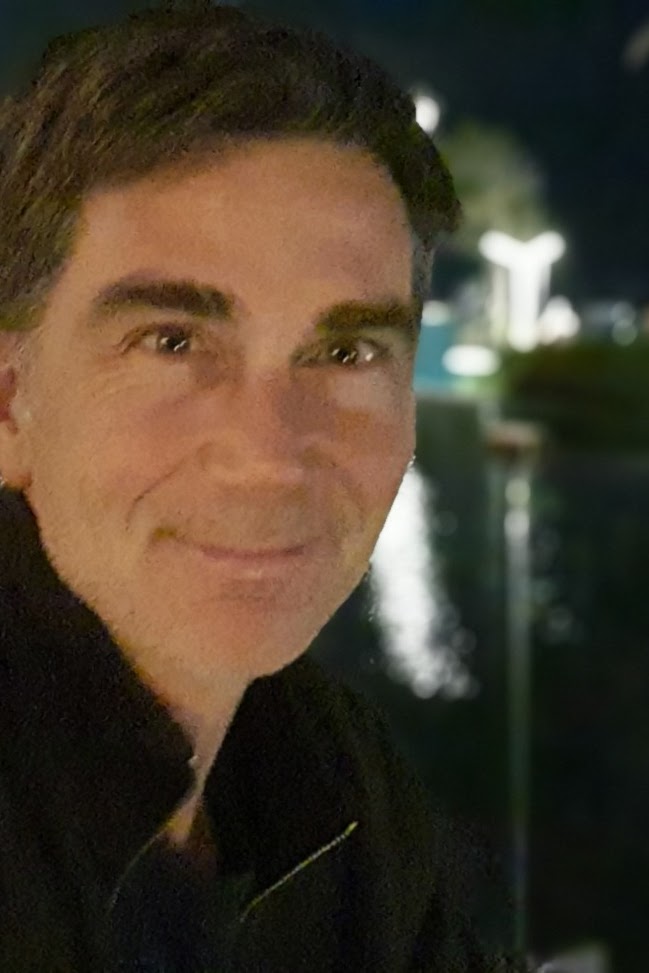}}]{Samson Lasaulce}
Samson Lasaulce is currently a Chief Research Scientist with Khalifa University, Abu Dhabi, UAE. He is the holder of the KU – TII 6G Chair on Native AI. He is also a Director of Research with CNRS, France. He has been a professor with the Department of Physics at École Polytechnique, France. His current research interests lie in distributed networks with a focus on game theory, optimal control, distributed optimization, machine learning for communication networks, energy networks, social networks, and decarbonization.
\end{IEEEbiography}

\begin{IEEEbiography}[{\includegraphics[width=1in,height=1.25in,clip,keepaspectratio]{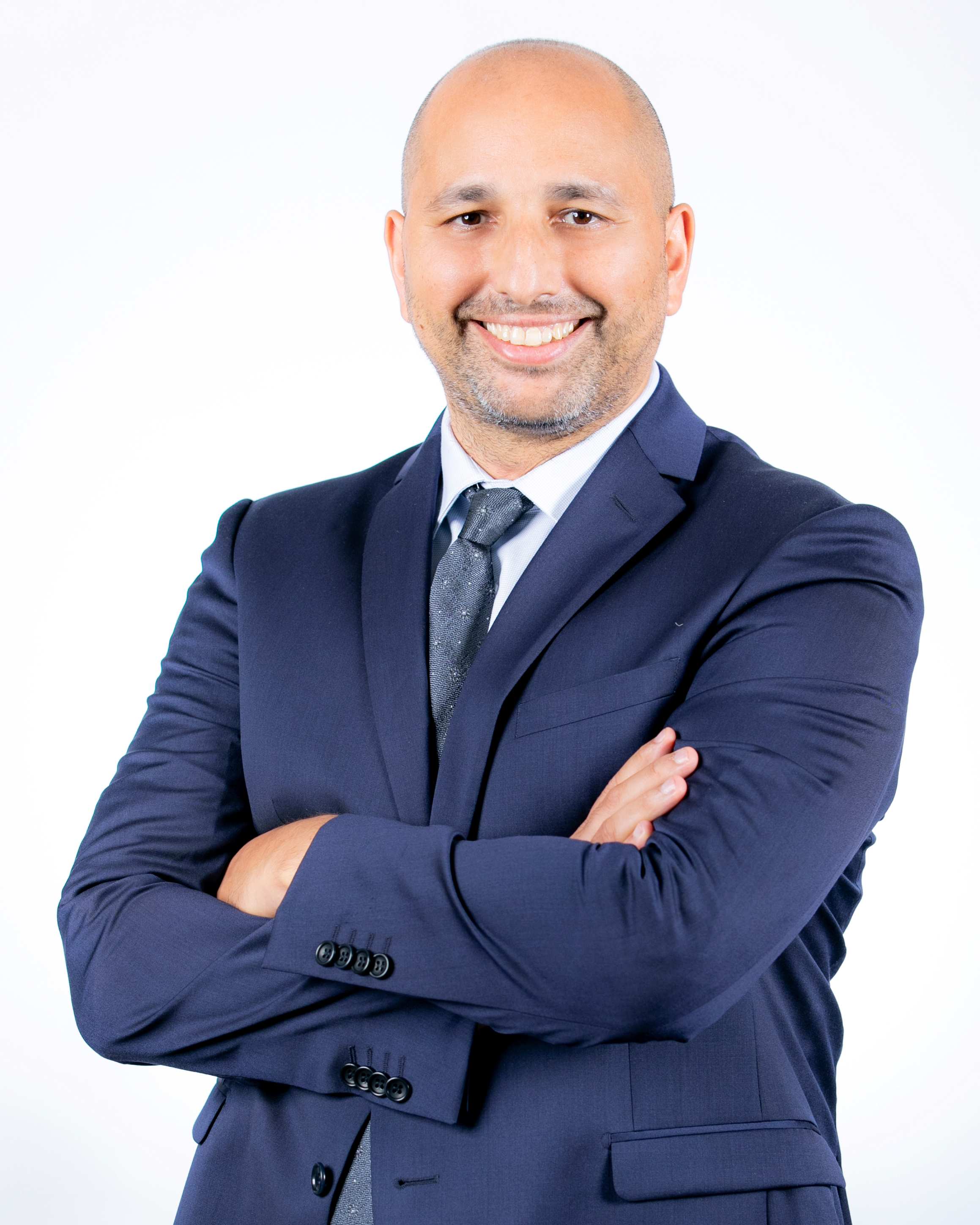}}]{Merouane Debbah}
	Mérouane Debbah is Professor at Khalifa University of Science and Technology in Abu Dhabi. He received the M.Sc. and Ph.D. degrees from the Ecole Normale Supérieure Paris-Saclay, France. Since 2023, he is a  Professor at  Khalifa University of Science and Technology in Abu Dhabi and founding director of the 6G center. He has managed 8 EU projects and more than 24 national and international projects. His research interests lie in fundamental mathematics, algorithms, statistics, information, and communication sciences research.
\end{IEEEbiography}

	\newpage

	\begin{figure*}
		\begin{center}
			\begin{large}
				{\bf Supplementary materials:} reproducing precoding matrices for the toy channel example of Figure \ref{fig:toy} in the main paper '{\it Per-antenna power constraints: constructing Pareto-optimal precoders with cubic complexity under non-negligible noise conditions}'
			\end{large}
		\end{center}
	\end{figure*}
	\begin{table*}[!th]
		\centering
		\caption{\centering Parameters of the highest-throughput (mean/sum) per-algorithm points on the toy real channel $H \in \mathbb{R}^{8 \times 3}$ given by Figure \ref{subfig:toy_ch_zf_slnr} - points marked with an arrow. {\bf (reference algorithms)}}
		\label{tab:throughput_reference}
		\begin{tabular}{|p{2.65cm}|p{4.5cm}|p{4.5cm}|}
			\hline
			{\bf \large Precoder type} & {\bf \large Zero-Forcing} & {\bf \large SLNR} \\
			\hline
			{\bf Best precoder $P$} 
			&
			$\begin{bmatrix}
				-0.4510 & -0.2924 & 0.2970 \\
				-0.4323 & -0.0362 & 0.9010 \\
				-0.0451 & 0.0812 & -0.1756 \\
				0.4504 & -0.0824 & -0.3683 \\
				-0.0111 & -0.0544 & -0.0910 \\
				0.0397 & -0.3846 & 0.6434 \\
				0.2615 & 0.7373 & -0.6228 \\
				-0.9188 & -0.3902 & -0.0556 
			\end{bmatrix}$
			&  
			$\begin{bmatrix}
				-0.3617 & -0.2835 & 0.1793 \\
				-0.3223 & 0.1973 & 0.9258 \\
				-0.0831 & 0.0911 & -0.1794 \\
				0.4284 & -0.2688 & -0.3639 \\
				-0.0111 & -0.0937 & -0.1181 \\
				0.1949 & -0.4364 & 0.6136 \\
				0.0409 & 0.8941 & -0.4457 \\
				-0.8573 & -0.4090 & -0.3118 
			\end{bmatrix}$
			\\
			\hline
			{\bf Weights \mbox{(Power allocation)}}
			&
			\mbox{$\kappa_1 = 0.3481$}
			\newline
			\mbox{$\kappa_2 = 0.2184$}
			\newline
			\mbox{$\kappa_3 = 0.4335$}
			&
			\mbox{$\kappa_1 = 0.2787$}
			\newline
			\mbox{$\kappa_2 = 0.3172$}
			\newline
			\mbox{$\kappa_3 = 0.4042$}
			\\
			\hline 
			{\bf Product} \newline $G = H^*P$ 
			& 
			$\begin{bmatrix}
				1.6993 & 0.0 & 0.0 \\
				0.0 & 1.3440 & 0.0 \\
				0.0 & 0.0 & 1.8115 
			\end{bmatrix}$
			& 
			$\begin{bmatrix}
				1.7195 & -0.3005 & 0.2597 \\
				-0.2720 & 1.9319 & 0.3773 \\
				0.2166 & 0.3476 & 2.0147 
			\end{bmatrix}$
			\\
			\hline
			{\bf SINR values}
			&
			\mbox{$\mathrm{SINR}_1 = 2.8878$}
			\newline
			\mbox{$\mathrm{SINR}_2 = 1.8063$}
			\newline
			\mbox{$\mathrm{SINR}_3 = 3.2814$}
			&
			\mbox{$\mathrm{SINR}_1 = 2.5537$}
			\newline
			\mbox{$\mathrm{SINR}_2 = 3.0683$}
			\newline
			\mbox{$\mathrm{SINR}_3 = 3.4758$}
			\\
			\hline
			{\bf Per-user throughputs}
			\newline
			\mbox{$\psi_k := $}
			\newline
			\mbox{$10 \log_{10}(1 + SINR_k)$}
			&
			\mbox{$\psi_1 = 5.8970$}
			\newline
			\mbox{$\psi_2 = 4.4814$}
			\newline
			\mbox{$\psi_3 = 6.3159$}
			\newline
			\bf{\large \mbox{$\psi_{avg} = 5.5647$}}
			&
			\mbox{$\psi_1 = 5.5069$}
			\newline
			\mbox{$\psi_2 = 6.0941$}
			\newline
			\mbox{$\psi_3 = 6.5087$}
			\newline
			\mbox{\large $\psi_{avg} = 6.0375$}
			\\
			\hline
		\end{tabular}
	\end{table*}
	\begin{table*}[!th]
		\centering
		\caption{\centering Parameters of the highest-throughput (mean/sum) per-algorithm points on the toy real channel $H \in \mathbb{R}^{8 \times 3}$ given by Figure \ref{subfig:toy_ch_zf_slnr}. The parameters correspond to feasible points marked with an arrow on Figures \ref{subfig:toy_p0}, \ref{subfig:toy_p1}, \ref{subfig:toy_p5}.
			\newline {\bf (proposed algorithm)}}
		\label{tab:throughput_algorithm}
		\begin{tabular}{|p{2.65cm}|p{4.5cm}|p{4.5cm}|p{4.5cm}|}
			\hline
			{\bf \large Iterations} & {\bf \large Zero (equal $\mu$)} & {\bf \large One $\mu$ update} & {\bf \large Convergence, $\delta = 0.01$} \\
			\hline
			{\bf Best precoder $P$} 
			&
			$\begin{bmatrix}
				-0.4077 & -0.3012 & 0.2314 \\
				-0.3905 & 0.0540 & 0.9188 \\
				-0.0659 & 0.0917 & -0.1773 \\
				0.4523 & -0.1650 & -0.3694 \\
				-0.0090 & -0.0723 & -0.1068 \\
				0.1244 & -0.4325 & 0.6251 \\
				0.1387 & 0.8426 & -0.5203 \\
				-0.8919 & -0.4063 & -0.1986 
			\end{bmatrix}$
			&  
			$\begin{bmatrix}
				-0.6967 & -0.5470 & 0.4171 \\
				-0.3274 & 0.0590 & 0.9086 \\
				-0.2792 & 0.4190 & -0.8600 \\
				0.6721 & -0.2845 & -0.6139 \\
				-0.0861 & -0.5395 & -0.8375 \\
				0.1426 & -0.5389 & 0.8237 \\
				0.1368 & 0.8337 & -0.5348 \\
				-0.8716 & -0.4168 & -0.2581 
			\end{bmatrix}$
			&  
			$\begin{bmatrix}
				-0.7252 & -0.6030 & 0.3279 \\
				-0.3299 & 0.0928 & 0.9308 \\
				-0.3893 & 0.4470 & -0.8046 \\
				0.7949 & -0.3489 & -0.4779 \\
				-0.3349 & -0.7205 & -0.6008 \\
				0.1922 & -0.6013 & 0.7750 \\
				0.2011 & 0.8479 & -0.4905 \\
				-0.8732 & -0.4647 & -0.1410 
			\end{bmatrix}$
			\\
			\hline
			{\bf Weights \mbox{(Lagrange, $\lambda_k$)}}
			&
			\mbox{$\lambda_1 = 0.3123$}
			\newline
			\mbox{$\lambda_2 = 0.2616$}
			\newline
			\mbox{$\lambda_3 = 0.4261$}
			&
			\mbox{$\lambda_1 = 0.2693$}
			\newline
			\mbox{$\lambda_2 = 0.2495$}
			\newline
			\mbox{$\lambda_3 =  0.4812$}
			&
			\mbox{$\lambda_1 = 0.3307$}
			\newline
			\mbox{$\lambda_2 = 0.3326$}
			\newline
			\mbox{$\lambda_3 = 0.3368$}
			\\
			\hline 
			{\bf Product} \newline $G = H^*P$ 
			& 
			$\begin{bmatrix}
				1.7385 & -0.1440 & 0.1383 \\
				-0.1755 & 1.6543 & 0.2224 \\
				0.0946 & 0.1247 & 1.9301 
			\end{bmatrix}$
			& 
			$\begin{bmatrix}
				2.0128 & -0.2445 & 0.2298 \\
				-0.2622 & 1.9765 & 0.3566 \\
				0.1133 & 0.1638 & 2.4909 
			\end{bmatrix}$
			& 
			$\begin{bmatrix}
				2.1391 & -0.2719 & 0.1532 \\
				-0.2653 & 2.1760 & 0.2744 \\
				0.1453 & 0.2666 & 2.2638 
			\end{bmatrix}$
			\\
			\hline
			{\bf SINR values}
			&
			\mbox{$\mathrm{SINR}_1 = 2.9065$}
			\newline
			\mbox{$\mathrm{SINR}_2 = 2.5335$}
			\newline
			\mbox{$\mathrm{SINR}_3 = 3.6363$}
			&
			\mbox{$\mathrm{SINR}_1 = 3.6413$}
			\newline
			\mbox{$\mathrm{SINR}_2 = 3.2667$}
			\newline
			\mbox{$\mathrm{SINR}_3 = 5.9677$}
			&
			\mbox{$\mathrm{SINR}_1 = 4.1696$}
			\newline
			\mbox{$\mathrm{SINR}_2 = 4.1328$}
			\newline
			\mbox{$\mathrm{SINR}_3 = 4.6920$}
			\\
			\hline
			{\bf Per-user throughputs}
			\newline
			\mbox{$\psi_k := $}
			\newline
			\mbox{$10 \log_{10}(1 + SINR_k)$}
			&
			\mbox{$\psi_1 = 5.9179$}
			\newline
			\mbox{$\psi_2 = 5.4821$}
			\newline
			\mbox{$\psi_3 = 6.6617$}
			\newline
			\bf{\large \mbox{$\psi_{avg} = 6.0206$}}
			&
			\mbox{$\psi_1 = 6.6664$}
			\newline
			\mbox{$\psi_2 = 6.3009$}
			\newline
			\mbox{$\psi_3 = 8.4309$}
			\newline
			\mbox{\large $\psi_{avg} = 7.1327$}	 
			&
			\mbox{$\psi_1 = 7.1346$}
			\newline
			\mbox{$\psi_2 = 7.1036$}
			\newline
			\mbox{$\psi_3 = 7.5527$}
			\newline
			\mbox{\large $\psi_{avg} = 7.2636$}
			\\
			\hline
		\end{tabular}
	\end{table*}	
	
	\begin{table*}[!th]
		\centering
		\caption{\centering High-precision values for a precoding matrix $P_{example}$ that is Pareto-optimal for channel from Figure \ref{subfig:toy_ch_zf_slnr} with $\omega_k = 0.04, \beta = 1$, but doesn't utilize full power at the last row.}
		\label{tab:example}
		\begin{tabular}{cc}
			
			$P_{example} =
			\begin{bmatrix}
				-0.131822459705 & -0.116783110219 & 0.984370125389 \\
				-0.179275984384 & -0.026748610369 & 0.983435118900 \\
				0.103075169577 & 0.169535064180 & -0.980119059316 \\
				0.228024544866 & -0.069956868783 & -0.971138941162 \\
				0.215688858238 & -0.045050174929 & -0.975422368191 \\
				-0.098145309120 & -0.178451934724 & 0.979041574716 \\
				0.059226063982 & 0.247555764804 & -0.967061743834 \\
				-0.572228658427 & -0.192375199215 & -0.257327740732 
			\end{bmatrix}$
			&
			$\begin{aligned}
				H^*P_{example} & = 
				\begin{bmatrix}
					0.86753830 & -0.00070736 & 0.01214739 \\
					-0.00119839 & 0.54445487 & 0.01513752 \\
					0.00007688 & 0.00005655 & 2.94679321 
				\end{bmatrix} \\
				\mathrm{SINR}_1 & \approx 4.305 \times 10^2, \\
				\mathrm{SINR}_2 & \approx 1.619 \times 10^2, \\
				\mathrm{SINR}_3 & \approx 5.427 \times 10^3.
			\end{aligned}$
		\end{tabular}
	\end{table*}

\end{document}